



\font\bm=msym10

\def\m#1{\hbox{\bm #1}}




\vsize=24.5truecm
\hsize=6.5truein


\nopagenumbers
\headline={\ifnum\pageno=1\hfil\else\hfil\rm\folio\fi}

\def\t#1{\widetilde{#1}}
\def\h#1{\widehat{#1}}
\def\b#1{\overline{#1}}
\def\oa#1{\buildrel #1\over\rightarrow}
\def\blowup#1{\,\t#1\oa\pi #1\,}
\def\norm#1{|\kern-1.72pt{|}\,#1\,|\kern-1.72pt{|}}
\def\sqr#1#2{{\vbox{\hrule height.#2pt \hbox{\vrule width.#2pt
  height#1pt \kern#1pt \vrule width.#2pt}\hrule height.#2pt}}}
\def\qed{\qquad \sqr 7 3}
\def\lcorner{\vrule height .4pt depth 0pt width 4pt \vrule height 6 pt depth
0pt width .4pt \;}
\def\bmatrix#1{\left[\matrix{#1}\right]}

\def\d{\partial}
\def\db{\bar\partial}
\def\dbd{\bar\partial\partial}
\def\ddb{\partial\bar\partial}

\def\lap{\triangle}
\def\tr{{\rm \,tr\,}}


\centerline{\bf Sequences of stable bundles over compact complex
surfaces.}

\bigskip
\centerline{N.~P.~Buchdahl}
\centerline{Department of Pure Mathematics}
\centerline{University of Adelaide}
\centerline{Adelaide, Australia 5005}
\centerline{e-mail: nbuchdah@maths.adelaide.edu.au} \vskip.75truein

\noindent{\bf 0. \quad Introduction.}
\bigskip

In [D2], Donaldson proved that every stable holomorphic bundle on
an algebraic surface equipped with a Hodge metric admits an
irreducible Hermitian-Einstein connection (unique up to gauge
equivalence).  The fact that bundles admitting Hermitian-Einstein
connections are quasi-stable (i.e., direct sums of stable bundles
all of the same slope) had already been proved by Kobayashi [Kob]
and L\"ubke [L]. Donaldson's result was subsequently generalized
by a number of authors, notably Uhlenbeck--Yau [UY] who extended
it to K\"ahler manifolds of arbitrary dimension, by Li--Yau [LY]
who extended it to general manifolds equipped with Gauduchon
metrics,  and by the author of this paper [B3], who proved the
Li--Yau result independently in complex dimension 2 by different
methods.

The identification of moduli spaces of stable bundles with moduli
spaces of anti-self-dual connections on an algebraic surface has
led to a number of extraordinary results in the differential
topology of 4-dimensional manifolds;  see, e.g., [D4], [FM1],
[FM2], [Kot], [OV]. A key theme in all these papers is that on an
algebraic surface, the Donaldson polynomials [D5] can be
expressed in terms of algebro-geometric quantities, thereby
facilitating their computation.   A frequently occurring problem
is that the moduli spaces involved are not compact, and although
there is a gauge theoretic compactification [D3] using the
results of Uhlenbeck [U1], [U2], this is not ostensibly an
algebraic object (but sometimes can be: see [Ma]). Instead,
authors have preferred to use Gieseker's compactification [Gie]
of  moduli spaces of stable bundles, but this also leads to some
difficulties in that Gieseker stability is not the same as
Mumford-Takemoto stability [OSS] and it is the latter definition
which is used in gauge theory. (Morgan [Mor] has however proved
that for the purposes of computing Donaldson polynomials, the
Gieseker compactification is sufficient.)

\medskip

The main result (Theorems 1.3 and 1.4) proved in this paper  is
that sequences of stable bundles of fixed topological type and
bounded degree have, after blowing-up and pulling-back
sufficiently often, strongly convergent subsequences.  Using
this, a natural topology can be defined on the space of stable
holomorphic structures of fixed degree on a given unitary bundle
and its pull-backs to blowups of the surface so that (under some
conditions) the space is compact---these results are given in the
related paper [B5].

The results presented here give an interpretation of the theorems
of Uhlenbeck in terms of the well-known phenomenon of jumping of
holomorphic structures.  In addition, the proof of Theorem~1.3
has a number of other useful spin-offs:  the main ones are the
derivation of quite precise and explicit estimates on the
curvature and connection forms of a concentrated anti-self-dual
connection in a ball (anti-self-dual with respect to a flat
metric),  and  monad descriptions of such connections. Indeed,
because much of the proof is local in nature, it comes comes
close to providing a compactification for moduli of
anti-self-dual connections on an arbitrary smooth 4-manifold $X$
in terms of strongly converging sequences of connections on
$X\#n\,\b{\m{CP}}_2$.

\bigskip
\medskip
The author is grateful to the Isaac Newton
Institute for Mathematical Sciences for its hospitality during
December 1994, and to the Department of Mathematics of the
University of Nantes for its hospitality during the first half of
1995 where this paper was completed.

\bigskip
\bigskip
\bigskip

\noindent{\bf 1. \quad Preliminaries.}

\bigskip
Let $X$ be a compact complex surface and let $\omega$ be a $\bar
\partial\partial$-closed positive (1,1)-form on $X$: it is a
theorem of Gauduchon [Gau] that every positive $(1,1)$-form has a
unique positive conformal rescaling such the rescaled form is
$\dbd$-closed and gives the same volume $V := Vol(X,\omega) :=
{1\over 2}\int_X \omega^2$.  With such a form $\omega$, the
degree  $deg(L) = deg(L,\omega)$ of a holomorphic line bundle $L$
on $X$ is unambiguously defined by the formula
$$
deg(L) \; := \;
{i\over 2\pi} \int_X f_L \wedge \omega\;,
$$
where $f_L$ is the
curvature of any hermitian connection on $L$.   The degree
depends only on $c_1(L)$ if and only if $b_1(X)$ is  even, and
when this is the case $\omega$ is cohomologous modulo the  image
of $\partial + \bar\partial$ to a closed form which itself is
unique up to the image of $\dbd$; ([B3], Proposition~2);   by the
Kodaira-Enriques classification of compact complex surfaces,
$b_1(X)$ is even iff $X$ admits a K\"ahler metric.

If $E$ is a holomorphic $r$-bundle on $X$, set $deg(E) :=
deg(\det\,E)$ and $\mu (E) := deg(E)/r$; the latter is called the
normalized degree  or {\it slope\/} of $E$.  A hermitian
connection on $E$ is {\it  Hermitian-Einstein\/} if the curvature
$F$ satisfies $\h F =  i\lambda \bf 1$ where $\h F := *\,(\omega
\wedge F) =: \Lambda F$,  $\lambda = (-2\pi / V )\!\cdot\! \mu
(E)$ and $\bf 1$ is the identity endomorphism of $E$.  The bundle
$E$ is ({\it semi-}) {\it stable\/} if $\mu ({\cal S}) <
(\le)~\mu (E)$ for every coherent subsheaf ${\cal S} \subset E$
with $0 < rank({\cal S}) < r$.  As mentioned in the introduction,
the main result of [B3] is that a bundle admits an irreducible
Hermitian-Einstein connection if and only if it is stable, this
generalizing the same result proved by Donaldson [D2] in the case
that $(X,\omega)$ is algebraic.  A  bundle admitting a
Hermitian-Einstein connection is a direct sum of  stable bundles
all of the same slope; i.e., is quasi-stable.

If $E$ has a Hermitian-Einstein connection with curvature $F$,
the equation $\omega \wedge \big( F - {1\over r}tr F\, {\bf 1}
\big) = 0$ and the skew-Hermitian property of $F$ give  $tr\big(
F - {1\over r}tr F \,{\bf 1} \big)^2 = | F - {1\over r}tr  F \,
{\bf 1}|^2\,dV$. Since the former $4$-form is a  representative
for the characteristic class  $8\pi^2\big( c_2 - {r-1\over
2r}c_1^2\big)(E)$, this motivates defining the {\it charge\/} of
$E $, $C(E)$, for an arbitrary $r$-bundle $E$ by the formula
$$
C(E) := \big( c_2 - {r-1\over 2r}c_1^2\big)(E) \;=\; {1\over
8\pi^2}\int_X tr\big( F - {1\over r}tr\, F \,{\bf 1} \big)^2
\quad  .\eqno (1.1)
$$
The charge is non-negative for any bundle
admitting a  Hermitian-Einstein connection, and when this is the
case, is  identically zero only if the induced Hermitian-Einstein
connection  on the adjoint bundle is flat.  Note that the charge
is  invariant under tensoring by line bundles: $C(E\otimes L) =
C(E)$  for any such $ L $.  In general, $C(E\otimes A) = aC(E) +
rC(A)$,  where $a,r$ are the ranks of $A,E$ respectively.

\medskip
As in [B3], it is often convenient to identify
holomorphic bundles with integrable unitary connections on a
fixed topological bundle.  (If $E_{\rm top}$ is such a bundle, a
unitary connection $A$ on $E_{\rm top}$ is integrable if the
$(0,2)$-component of the curvature $F(A)$ vanishes.)  The group
of complex automorphisms of $E_{\rm top}$ acts on the set of
unitary connections on $E_{\rm top}$ via
$$
(g,A) \mapsto g\cdot
A \qquad {\rm where} \quad d_{g\cdot A} :=  g\circ \db_A \circ
g^{-1} + {g^*}^{-1}\circ \d_A \circ  {g^*}\;,
$$
and the
corresponding curvatures are related by
$$
F^{0,2}(g\cdot
A)=gF^{0,2}(A)g^{-1}, \qquad F^{1,1}(g\cdot
A)=g[F(A)+\db_A(h^{-1}\d_A h)]g^{-1}\qquad {\rm for }\quad h :=
g^*g\,.
$$
The action preserves the set of integrable
connections, and two such connections induce isomorphic
holomorphic structures iff they lie in the same orbit.

\medskip
The metrics on blowups of $X$ which will be used to
define stability are as described in [B3],  recalled here briefly
for convenience. For more details, see [B5,\S3]. In the
discussion which follows, it is  not assumed that $X$ is compact.

\smallskip

Let $\blowup X$ be the blowup of $X$ at $x_0 \in X$.  Let $L_0 :=
\pi^{-1}(x_0)$ be the exceptional divisor so $\pi^*\omega$ is
everywhere non-negative and is degenerate only in directions
tangent to $L_0$.  Let $\sigma$ be $i/2\pi$ times the curvature
form  of any hermitian connection on the line bundle ${\cal
O}(-L_0) =:  {\cal O}(1)$ restricting positively to $L_0$, and
let $\omega_{\epsilon} := \pi^*\omega + \epsilon \sigma$ for
$\epsilon  > 0$; (recall ${\cal O}(L_0)\!\mid_{L_0} \simeq {\cal
O}_{L_0}(-1)$).  It  follows that if $\epsilon$ is sufficiently
small then $\omega_{\epsilon}$ defines a positive form in a
neighbourhood of  $L$; if $\sigma$ is compactly supported in $\t
X$ then  $\omega_{\epsilon}$ is everywhere positive for
sufficiently small  $\epsilon$. Note that $\sigma$ can be taken
to have support in neighbourhoods of $L_0$ of the form
$\pi^{-1}(U)$ where $U$ is an arbitrarily small neighbourhood of
$x_0$.  If $\omega$ is $\dbd$-closed and $\sigma$ is compactly
supported, it follows from the fact that $L$ has
self-intersection  $-1$ that $Vol(\t X,\omega_{\epsilon}) =
Vol(X,\omega) -\epsilon^2/2$, and if $\omega$ is $d$-closed, then
so too is $\omega_{\epsilon}$.

Now let $\blowup X$ be a modification of $X$ consisting of $n$
blowups, and let $\sigma_i$ be a closed smooth $(1,1)$-form on
$\t X$ corresponding as in the last paragraph to the $i$-th
blowup (pulled-back to $\t X$).  Let $\m R^n_+ :=\{ \alpha =
(\alpha_1,\ldots,\alpha_n) \in {\m{R}}^n\!\mid \alpha_i >0 \,,\;
i  = 1,\ldots,n \}$ and for $\alpha \in \m R^n_+\,$ let
$\rho_{\alpha} := \sum \alpha_i\sigma_i$, so
$\rho_{\alpha}\!\cdot\!\rho_{\alpha} = -\sum\alpha_i^2 =:
-|\alpha|^2$; (this definition differs slightly  from that in
[B3] where $\rho$ has the opposite sign).  The vector $\alpha \in
\m R^n_+$ will be called {\it suitable\/} if $\omega_{\alpha} :=
\pi^*\omega+\rho_{\alpha}$ defines a positive $\dbd$-closed
$(1,1)$-form on $\t X$ for some choice of the forms $\sigma_i$.
If $\t X'$ is another modification of $X$ consisting of $n$
blowups, each $\sigma_i$ on $\t X$ corresponds to a $\sigma_i'$
on $\t X'$ in $H^2(\t X,\m R)$ since $\t X$ and $\t X'$ are
diffeomorphic so the  degrees of line bundles on $\t X$ and $\t
X'$ defined by $\pi^*\omega+\rho_{\alpha}$ and
$\pi'^*\omega+\rho'_{\alpha}$ are the same.   The notation will
sometimes be abused by using the one symbol $\omega_{\alpha}$ to
denote a metric on a blowup, even though the blowup may be
changing.

\medskip

The following result summarises some of the relationships between
stability on $X$ and $\t X$;  the proof is given in [B5,\S3]:

\bigskip
\noindent {\bf Proposition~1.2.} \quad {\sl
\smallskip
\itemitem{(a)~~} If $\t E$ is a bundle on $\t X$ such that
$\pi_*\t E$ is stable, then there exists $\epsilon>0$ such that
$\t E$ is $\omega_{\alpha}$-stable for all suitable $\alpha\in \m
R^n_+$ satisfying $|\alpha|<\epsilon$;
\itemitem{(b)~~} If $\t E$
on $\t X$ is $\omega_{\epsilon\alpha}$-semi-stable for all
$\epsilon>0$ sufficiently small, then $\pi_*\t E$ is semi-stable;
\itemitem{(c)~~} If $b_1(X)$ is even, then for any $r_0, C_0$
there exists $\epsilon_0=\epsilon_0(\omega,r_0,C_0)>0$ with the
property that any bundle $\t E$ on $\t X$ of rank $\le r_0$ and
charge $\le C_0$ which is $\omega_{\alpha}$-stable for some
suitable $\alpha \in \m R^n_+$ satisfying $|\alpha|<\epsilon_0$
is also $\omega_{\epsilon\alpha}$-stable for all $\epsilon\in
(0,\epsilon_0/|\alpha|)$.  \quad  \qed \par}

\bigskip
\noindent In particular, the pull-backs of stable
bundles are stable, and when $b_1(X)$ is even, the pull-backs of
stable bundles all of the same topological type are stable with
respect to the same metric on the blowup.

Note that for any bundle $\t E$ on $\t X$, the direct image
$\pi_*\t E$ is a torsion-free sheaf on $X$.  The double-dual
$(\pi_*\t E)^{**}$ is reflexive, hence locally free [OSS], and a
torsion-free sheaf is (semi-)stable iff its double-dual is.  It
is not hard to show that $C((\pi_*\t E)^{**}) \le C(\t E)$ with
equality iff $\t E=(\pi^*E)\otimes L$ for some bundle $E$ on $X$
and some line bundle $L$ on $\t X$ trivial off the exceptional
divisor; ([B5, Proposition~2.9]).

\bigskip

Let $\{E_i\}$ be a sequence of stable bundles on $(X,\omega)$ all
of the same topological type and of uniformly bounded degree.
Then by Uhlenbeck's weak compactness theorem [U1], [S],  there is
a finite set $S\subset X$ such that, after gauge transformations,
a subsequence of the corresponding sequence of Hermitian-Einstein
connections converges weakly in $L^p_2$ for some fixed $p>4$ and
strongly in $C^1$  on compact subsets of $X\backslash S$ to a
Hermitian-Einstein connection which, by removability of
singularities [U2], then extends across $S$ to define a
Hermitian-Einstein connection on a bundle $E$ with $C(E) \le
C(E_i)$, endowing $E$ with the structure of a quasi-stable
holomorphic bundle on $(X,\omega)$.  This type of convergence
will be summarised by saying that the subsequence of $\{E_i\}$
{\it converges weakly\/} to $E$ (off $S$) (with respect to
$\omega$).  If $S$ is empty, the convergence will be called
``strong convergence" (with respect to $\omega$).

Suppose for the moment that $b_1(X)$ is even and that the weak
limit $E$ is stable.  Blow up $X$ along $S$ to $\t X$ and fix a
metric on $\t X$ of the form $\omega_{\alpha}$, where $|\alpha| <
\epsilon_0$ as in Proposition~1.2(c).  Then $\pi^*E_i$ is
$\omega_{\alpha}$-stable for all $i$ so by weak compactness on
$(\t X,\omega_{\alpha})$, there is a finite set $\t S \subset \t
X$ such that a subsequence of $\{\pi^*E_i\}$ converges weakly to
some $\omega_{\alpha}$-quasi-stable bundle $\t E$ on $\t X$.
Using Lemmas 2.1 and 2.2 of the next section, it is not hard to
show that under these circumstances, $\t S \subset \pi^{-1}(S)$,
$\t E$ is stable, and $(\pi_*\t E)^{**}=E$.

If it could be shown that $\t E$ were non-trivial on the
exceptional divisor, then the amount of charge bubbled by the
sequence on the blowup would be strictly less than that bubbled
on $X$, and therefore by induction on this amount, iteration of
the procedure would eventually lead to a strongly convergent
sequence on some blowup.  Unfortunately, there is no guarantee
that the amount bubbled on the blowup is not the same as that
bubbled on the original surface, and indeed examples can be
constructed where this occurs. However, if instead of blowing up
$X$ along the fixed set $S$, the points in the center of the
modification are allowed to vary (but converge to $S$), then
strongly convergent subsequences can be obtained:

\bigskip
\noindent {\bf Theorem~1.3.} \quad {\sl Let $X$ be a
compact complex surface equipped with a positive $\dbd$-closed
$(1,1)$-form $\omega$.  Let $\{A_i\}$ be a sequence of
Hermitian-Einstein connections on a fixed unitary bundle $E_{\rm
top}$ of rank $r$ over $X$ such that the corresponding
holomorphic bundles $E_i$ are stable and are of uniformly bounded
degree.  Suppose that $E_i$ converges weakly to $E$ off $S
\subset X$.  Then there is a subsequence $\{E_{i_j}\}\subset
\{E_i\}$ such that:
\medskip
\item{1.\quad} There is a sequence
of blowups $\t X_{i_j}$ of $X$ consisting of  $n \le 2C(E_{\rm
top})-2C(E)-1$  blowups at simple points converging to a blowup
$\blowup X$ of $X$;
\item{2.\quad} The exceptional divisor of $\t
X$ is $\pi^{-1}(S)$;
\item{3.\quad} There are complex
automorphisms $g_{i_j}$ of $\pi_{i_j}^*E^{}_{\rm top}$ such that
$\{g_{i_j}^{}\cdot (\pi_{i_j}^*A_{i_j}^{})\}$ converges strongly
in $C^1$ to define a holomorphic bundle $\t E$ on $\t X$ with
$(\pi_*\t E)^{**} = E$;
\item{4.\quad} For each compact set $K
\subset X\backslash S$ there is a constant $C_K$ such that
$\sup_K(|g_{i_j}^{}|+|g_{i_j}^{-1}|) \le C_K$ for all $i_j$, with
$\{g_{i_j}\}$ converging uniformly in $C^2(K)$;
\item{5.\quad}
For any choice $T$ of $r^2+1$ disjoint balls in $X\backslash S$
and sufficiently small suitable $\alpha \in \m R^n_+$ the
connections $g_{i_j}^{}\!\cdot (\pi_{i_j}^*A_{i_j}^{})$ can be
taken to be Hermitian-Einstein with respect to
$\pi_{i_j}^*\omega+\rho_{\alpha}^{}$ on the complement of $T$;
\item{6.\quad} If $E$ is stable, then for any suitable $\alpha
\in \m R^n_+$ sufficiently small, the bundles
$\pi_{i_j}^*E_{i_j}$ are
$(\pi_{i_j}^*\omega+\rho_{\alpha}^{})$-stable, $\t E$ is
$\omega_{\alpha}$-stable, and the connections $g_{i_j}^{}\!\cdot
(\pi_{i_j}^*A_{i_j}^{})$ can be taken to be Hermitian-Einstein
with respect to $\pi_{i_j}^*\omega+\rho_{\alpha}^{}$. \par }

\smallskip

\noindent The convergence of a sequence of blowups should be
interpreted as the convergence of a sequence of integrable
complex structures on the  same underlying smooth manifold $X\#n
\b{\m P}_2$ endowed with a fixed Riemannian metric.  In the case
in  hand, these complex structures are isomorphic on the
complement of  a fixed open set with strictly pseudoconvex
boundary.

\bigskip
It is evident from the statement of the theorem that
complications arise when the weak limit $E$ is only semi-stable
but not stable.  However, for bundles of rank $2$,  these
difficulties can be overcome:

\bigskip
\noindent {\bf Theorem~1.4.} \quad {\sl Under the
hypotheses of Theorem~1.3, suppose in addition that $b_1(X)$ is
even and that $E_{\rm top}$ has rank $2$.  Then there is a
subsequence $\{E_{i_j}\}\subset \{E_i\}$ such that:

\medskip
\item{1.\quad} There is a sequence of blowups $\t X_{i_j}$ of $X$
consisting of  $n \le 2C(E_{\rm top})-1$  blowups at simple
points converging to a blowup $\blowup X$ of $X$;
\item{2.\quad}
For some suitable $\alpha \in \m R^n_+$, $\pi_{i_j}^*E_{i_j}$ is
$(\pi_{i_j}^*\omega+\rho_{\alpha})$-stable and the corresponding
sequence of Hermitian-Einstein connections converges strongly on
$\t X$ to define an $\omega_{\alpha}$-stable bundle $\t E$ on $\t
X$;
\item{3.\quad} $\det \t E=\pi^*\det E$, $\pi_*\t E$ is
semi-stable, and there are non-zero homomorphisms $(\pi_*\t
E)^{**} \to E$,  $E \to (\pi_*\t E)^{**}$.\par}

\medskip
\noindent Note that the intertwining operators $g_{i_j}$
linking the two sequences of connections can no longer be
guaranteed to enjoy the same uniform bounds as given by the
fourth statement in Theorem~1.3.

\bigskip

Let $E_{\rm top}$ be a unitary $r$-bundle over $(X,\omega)$ and
let ${\cal M}(X,E_{\rm top})$ denote the space of isomorphism
classes of quasi-stable holomorphic structures on $E_{\rm top}$.
Consider the set of pairs $(\t X, \t E)$ where

\smallskip{\sl
\itemitem{{\rm 1. } } $\t X$ is a blowup of $X$;
\itemitem{{\rm 2. } } $\t E$ is a holomorphic bundle on $\t X$
topologically isormorphic to $\pi^*E_{\rm top}$ such that
$\pi_*E$ is semi-stable;
\itemitem{{\rm 3. } } $\t E$ is
$\omega_{\alpha}$-quasi-stable for all suitable $\alpha$ in an
open set of such;
\itemitem{{\rm 4. } } If $b_1(X)$ is odd, $\deg(\pi_*\t E,\omega)=0$.
\par}

\smallskip
\noindent  On this set, define an equivalence relation
$\sim$ by setting $(\t X,\t E) \sim (X',E')$ iff there is a joint
blowup $\t X'$ such that the pull-backs of $\t E$ and $E'$ to $\t
X'$ are isomorphic in the usual sense, and let  $\b{\cal
M}(X,E_{\rm top})$ denote the set of equivalence classes.

A topology is defined on $\b{\cal M}$ by defining $\{[(\t X_i,\t
E_i)]\}\subset \b{\cal M}$ to converge to $[(\t X,\t E)]$ iff
$[(\t X_i,\t E_i)]$ can be represented by  a sequence of blowups
$\t X_i$ converging to $\t X$ with $\omega_{\alpha}$-stable
bundles $\t E_i$ on $\t X_i$ converging strongly to $\t E$ on $\t
X$.

\smallskip
The subset ${\cal M}(X,E_{\rm top}) \subset \b{\cal
M}(X,E_{\rm top})$ consists of the bundles on blowups which are
holomorphically trivial on the exceptional divisor.  This is an
open set for if $\{\t E_i\} \subset \b{\cal M}\backslash {\cal
M}$ converges to $\t E$ then each $\t E_i$ in a subsequence is
non-trivial on some irreducible component $L$ of the exceptional
divisor and by semi-continuity of cohomology, the fact that
$H^0(L,\t E_{i_j}(-1))\not=0$ implies $H^0(L,\t E(-1))\not=0$.
Theorems 1.3 and 1.4 imply that under certain restrictions, every
sequence in ${\cal M}$ has a subsequence converging in $\b{\cal
M}$.

In [B5,\S6] it is proved that if $b_1(X)$ is even and
$rank(E_{\rm top})=2$ then $\b{\cal M}(X,E_{\rm top})$ is compact.
It is also proved that if every bundle $E$ on $X$ with
$rank(E)=rank(E_{\rm top})$, $c_1(E)=c_1(E_{\rm top})$ and $C(E)
\le C(E_{\rm top})$ which is semi-stable is actually stable, then
$\b{\cal M}(X,E_{\rm top})$ is  compact and Hausdorff space. If
$b_1(X)$ is even and $c_1(E_{\rm top})$ is not a torsion class in
$H^2(X,\m Z_r)$ for $r=rank(E_{\rm top})$, this condition will be
satisfied for generic $\omega$.  If $b_1(X)$ is odd, it will only
be satisfied if $c_2(E_{\rm top})<0$ (since $(L^*\otimes\det
E)\oplus {\cal O}^{r-1}$ is semi-stable if $c_1(L)=0$ and
$\deg(L)=\deg(E)/r$; the existence of such $L$ follows from the
proof of Proposition~2 in [B3]).

\bigskip
To outline the structure of the remainder of this paper,
let $(X,\omega)$ be a given compact hermitian surface with
$\dbd\omega=0$, and let $\{A_i\}$ be a sequence of
Hermitian-Einstein connections on a unitary $r$-bundle $E_{\rm
top}$ over $X$ which is converging in $C^1$ on compact subsets of
the complement of some non-empty finite set $S \subset X$.

Suppose $x_0\in S$. After twisting $A_i$ by a converging sequence
of Hermitian-Einstein connections on a trivial line bundle in a
neighbourhood of $x_0$, it can be assumed without loss of
generality that $A_i$ is anti-self-dual in this neighbourhood.
Assume for the moment that the metric $\omega$ is the standard
flat metric in this neighbourhood, and dilate to the ball
$B(x_0,5)$ of radius $5$. Transferring everything to $\m C^2$
gives  a sequence of anti-self-dual connections $A_i$ on an
$SU(r)$ bundle over the ball which are converging in $C^1$ on
compact subsets of $B(0,5)\backslash\{0\}$.  By dilating further
if necessary, it can be assumed that the limiting connection $A$,
after removal of singularities, satisfies $||F(A)||_{L^2}^2 \le
\epsilon$ where $\epsilon >0 $ can be prescribed arbitrarily
small, and $lim\big|\,||F(A_i)||^2-8\pi^2k\big| \le \epsilon$ for
some integer $k>0$.

In \S2, it will be shown that (perhaps after further dilations)
the connections $A_i$ can be truncated by a fixed cutoff function
$\rho$ and the resulting connections $\rho A_i$ can be perturbed
to anti-self-dual connections $\rho A_i+a_i$  on a bundle over
$S^4$ with $c_2=k$ (i.e., instantons).  Such connections have an
explicit description in terms of monads on $\m P_2$, that is, in
terms of linear algebraic data ([D1]), and the connections $\rho
A_i$ can be estimated by explicit calculation of the instanton
connections;  this is performed in \S\S3--4.

Although $\rho A_i+a_i$ provides a good approximation to $\rho
A_i$, it is unfortunately not the case that after pulling back to
the blowup of $\m C^2$ at the origin, the Hermitian-Einstein
connection corresponding to $\pi^*(\rho A_i+a_i)$ necessarily
provides a good approximatio to that for $\pi^*(\rho A_i)$:  the
reason for this is that although bounded functions pull back to
bounded functions on the blowup, $L^2$ functions do not pull-back
to $L^2$ functions.  Thus in \S5 a further perturbation is
required to obtain a sequence of instantons on $S^4$ which does
have this good approximation property, and this then reduces the
convergence problem for $\{\rho A_i\}$ to that of a sequence of
genuine instantons on $\m C^2$.

In \S6 it is shown that the main compactness result holds for
sequences of instantons on $S^4$ viewed as holomorphic bundles on
$\m P_2$ trivial on the line at infinity, with the condition of
stability suppressed by fixing trivialisations on this
line---this is a purely algebraic calculation involving monads on
$\m P_2$ and a blowup at a point.  Following this, the last
pieces of the proof of Theorem~1.3 are put in place.  The last
section is devoted to the proof of Theorem~1.4, which consists
largely of sheaf theory.

\bigskip
\bigskip
\bigskip

\noindent{\bf 2. \quad Transferring concentrated connections.}
\bigskip

Let $X$ be a compact complex surface, $\omega$ be a $\dbd$-closed
positive $(1,1)$-form on $X$, and let $\{A_i\}$ be a sequence of
irreducible Hermitian-Einstein connections on a unitary bundle
over $X$ which are converging in $C^1$ on compact subsets of
$X\backslash S$ for some finite set $S\subset X$.

In \S1, it was assumed that the metric $\omega$ was flat in a
neighbourhood of each of the points of $S$.  To see that this is
possible, let $x_0$ be a point of $S$ and choose holomorphic
coordinates $z=(z^0,z^1)$ in a neighbourhood of $x_0$ such that
$\omega$ agrees with the  standard metric $\omega_0 =
(i/2)\ddb|z|^2$ at $x_0$.  Since $\omega$ is $\dbd$-closed, there
is a $(1,0)$-form $u$ in a neighbourhood of $x_0$ such that
$\omega=\omega_0 + \db u + \d \bar u$ in that neighbourhood, and
a simple Taylor  expansion argument shows that $u$ can be chosen
to vanish to second  order at $x_0$.  Then if $\varphi \colon \,
[0,\infty) \to [0,1]$ is a smooth function which is $0$ on
$[0,r_1]$ and $1$ on $[2r_1,\infty)$, for sufficiently small
$r_1>0$ the form $\omega' := \omega_0 + \db[\varphi(|z|^2) u]+
\d[\varphi(|z|^2) \bar u]$ defines a $(1,1)$-form on $X$ which is
everywhere positive, is $\dbd$-closed, agrees with the standard
metric in a neighbourhood  of $x_0$ and with $\omega$ on the
complement of a slightly larger neighbourhood, and is homologous
to $\omega$ modulo the image of $\d +  \db$.  Hence the notions
of (semi-)stability for the two metrics coincide.

\smallskip
The following result shows that under certain
circumstances, the notions of weak convergence for the two
metrics also coincide---in this lemma, it is not necessary to
assume that $X$ is compact:

\bigskip

\noindent{\bf Lemma~2.1.}\quad  {\sl Let $\omega$, $\omega'$ be
positive $(1,1)$-forms on $X$.   Let $\{A_i\}$, $\{A_i'\}$ be
sequence of unitary connections on a bundle such that
$A_i'=g_i\cdot A_i$ for some complex automorphism, and suppose
that the following conditions hold:
\medskip
\itemitem{1.~~}
$||F(A_i)||_{L^2}+||F(A_i')||_{L^2}$ is uniformly bounded;
\itemitem{2.~~} $||\omega\wedge
F(A_i)||_{L^p}+||F^{0,2}(A_i)||_{L^p} +||\omega'\wedge
F(A_i')||_{L^p} +||F^{0,2}(A_i')||_{L^p}$ is uniformly bounded
for some $p>4$;
\itemitem{3.~~} There are finite subsets
$S,\,S'\subset X$ such that $\{A_i\}$ converges in $C^1$ on
compact subsets of $X\backslash S$, and $\{ A'_i\}$ converges in
$C^1$ on compact subsets of $X\backslash S'$;
\itemitem{4.~~}
After gauge transformations, the limits $A_{\infty}$,
$A'_{\infty}$ extend to connections on some other unitary bundles
over $X$. \par

\smallskip
\noindent If $\sup_Q(|g_i|+|g_i^{-1}|)$ is uniformly
bounded for some $Q$ which is the complement of a union of
non-intersecting closed balls containing the points of $S\cup S'$
then

\smallskip
\itemitem{1.~~} $S=S'$;
\itemitem{2.~~}  a subsequence
of $\{g_i\}$ converges in $C^2$ on compact subsets of
$X\backslash S$ to some $g_{\infty}$;
\itemitem{3.~~}  there are
unitary automorphisms $u,u'$ such that $ug_{\infty}u'^{-1}$
extends to a $C^2$  isomorphism of bundles over $X$;
\itemitem{4.~~} for any $\epsilon >0$ there exists $r>0$ such
that $ \big|\int_{B(x_0,r)}(|F(A_i)|^2-|F(A_i')|^2)\,dV\,\big|
\le \epsilon $ for all $i$ and every $x_0\in S$;
\itemitem{5.~~}
if $\omega=\omega'$ then $\sup_X(|g_i|+|g_i^{-1}|)$ is uniformly
bounded.\par}

\medskip
\noindent{\bf Proof:}\quad If $x_0\in S'\backslash S$,
then in some small ball around $x_0$ the sequence $\{A_i\}$ is
converging.  The uniform $L^p$ bound on $\omega' \wedge
F(g_i\cdot A_i)$ together with the Maximum Principle and the
bound on $|g_i|+|g_i^{-1}|$ on the boundary of this ball then
give a uniform $C^0$ bound on $|g_i|+|g_i^{-1}|$ over this ball.
By Lemma~2.1 of [B4] the hypotheses imply  that $x_0$ is not in
fact a bad point for the sequence $\{A_i'\}$ so $S' \subset S$,
and by symmetry it follows $S=S'$.

The uniform bounds on $g_i$ and $g_i^{-1}$ together with the
convergence of the sequence of connections implies that a
subsequence of $\{g_i\}$ converges in $C^2$ on compact subsets of
$X\backslash S$ and the existence of  $u, u'$ follows from the
existence the gauge transformations enabling the limiting
connections to extend, together with Hartogs' Theorem.

To prove the fourth statement, choose $r$ so small that
$|\int_{B(x_0,r)}(|F(A_{\infty})|^2+|F(A'_{\infty})|^2)\,dV|\le
\epsilon/3$.  By Stokes' Theorem~and the fact that on  the
boundary of $B(x_0,r)$ the connections are converging in $C^1$ it
follows that  $\big|\int_{B(x_0,r)}(\tr F(A_i)^2-\tr
F(A_i')^2)\big| \le \epsilon/2$ for $i$ sufficiently large.  With
respect to $\omega$,
$$
\eqalign{ ||F(A_i)||^2_{L^2(B(x_0,r))} &=
\int_{B(x_0,r)}\tr F(A_i)^2+  2|| F_+(A_i)||^2_{L^2(B(x_0,r))}
\cr
&\le \int_{B(x_0,r)}\tr F(A_i)^2 + Const. r^{4(p-2)/p}||
F_+(A_i)||^2_{L^p(B(x_0,r))}\;, }
$$
and because of the uniform
bound on the $L^p$ norm of $F_+(A_i)$, the second term can be
made arbitrarily small by choosing $r$ small enough.  The same
argument applies to the other sequence using $\omega'$, and the
result follows since the two metrics compare uniformly on $X$.

The last statement follows easily from the Maximum Principle, the
local existence of uniformly bounded solutions to $\lap u = |\h
F(A_i)|$, the Sobolev Embedding theorem and Theorem~9.20 of [GT].
\qquad \qed

\bigskip

To find intertwining operators $g_i$ satisfying the hypotheses of
Lemma~2.1 usually involves stability of either of the two limits,
together with the following version of semi-continuity of
cohomology:

\bigskip

\noindent {\bf Lemma~2.2.}\quad \sl Let $Q$ be an open subset of
$X$ with (possibly empty) smooth strictly pseudo-concave
boundary.  Let $\{A_i\}$ be a sequence of integrable connections
on a unitary bundle converging weakly in $L^p_2$ and uniformly in
$C^1$ on a neighbourhood of $\bar Q$ to $A_{\infty}$.  If $s_i$
are non-zero  sections over $Q$ satisfying $\db_i s_i = 0$ and
$|| s_i||_{L^2(Q)} = 1$, then there exists a subsequence
converging in $C^2$ on a neighbourhood of $\bar Q$ to a section
$s_{\infty}$ satisfying $\db_{\infty} s_{\infty} = 0$ and
$||s_{\infty}||_{L^2(Q)} = 1$.   \rm

\medskip
\noindent {\bf Proof: }  On any compact set $K \subset
Q$, ellipticity of the $\db$-operator and convergence of the
sequence of connections gives a uniform $C^0(K)$ bound on $s_i$
in terms of an $L^2(Q)$ bound. In a sufficiently small
neighbourhood of any $x_0\in \partial Q$ the sequence of sections
can be viewed as a sequence of vector-valued holomorphic
functions  over this neighbourhood.  Strict pseudoconvexity of
$\partial (X\backslash Q)$  and the fact that $dim(X) > 1$ imply
that these functions extend  across the boundary, and Cauchy
estimates give $C^0$ bounds on the extensions in terms of $C^0$
bounds on a compact subset of $Q$.  Hence there is an a  priori
$C^0$ bound on the sections $s_i$ over a fixed neighbourhood of
$\bar Q$.  Ellipticity of $\db$ and convergence of the  sequence
of connections then gives a uniform $L^p_3$ bound on $s_i$ over a
slightly smaller neighbourhood and hence there is a  subsequence
$\{s_{i_j}\}$ converging strongly in $C^2(\bar Q)$ using
compactness of the embedding $L^p_3 \subset C^2$ for $p>4$. \qed

\bigskip

Applying these results to the sequence $\{A_i\}$, let $E_i$ be
the holomorphic bundle defined by $A_i$ and suppose that the weak
limit $E$ defined by the sequence is also stable.  Since each
$E_i$ is also stable with respect to the flattened metric
$\omega'$, there is another finite set $S'\subset X$ such that
the sequence of $\omega'$-Hermitian-Einstein connections
$\{A'_i\}$ has a subsequence converging weakly off $S'$ to define
an $\omega'$-quasi-stable bundle $E'$.  Choose a disjoint union
of balls surrounding the points of $S\cup S'$ and let $Q$ be the
complement.  If $g_i$ is an intertwining operator (uniquely
determined up to scale) such that $A_i'=g_i\cdot A_i$,
renormalise $g_i$ so that $||g_i||_{L^2(Q)}=1$.  By Lemma~2.2 a
subsequence converges in $C^2$ to a non-zero limit $g_{\infty}$
which then extends to the whole of $X$ by Hartogs' Theorem.
Since $E$ is stable and $E'$ is semi-stable (with respect to
either metric), $g_{\infty}$ is an isomorphism, and therefore
Lemma~2.1 now applies to show that $S=S'$ and that the amount of
charge bubbled by each sequence is the same at each point of $S$.

\bigskip

Let $A$ be a typical connection in the sequence $\{A_i\}$ (where
it is now assumed that the metric $\omega$ is standard in a
neighbourhood of each point of $S$). After twisting the
connection by a fixed Hermitian-Einstein connection on a trivial
line bundle in a neighbourhood of $x_0$, it can be assumed
without loss of generality that $A$ is anti-self-dual; (if
$b_1(X)$ is odd, it should be assumed from the outside that every
connection in the sequence has degree $0$).  Fix a small number
$\eta>0$ and choose $r>0$  sufficiently small that the $L^2$ norm
of $F(A)$ over the annulus  $B(x_0,2r)\backslash B(x_0,r)$ is
less than $\eta$.  If $\eta$ is  sufficiently small, it follows
as in Chapter~9 of [FU] that after pulling  back the connection
from $B(x_0,4r)$ to $B(x_0,4)$ there is a gauge  transformation
over the annulus $B(x_0,4)\backslash B(x_0,3)$ so that the
gauged connection, also denoted by $A$, has its $C^1$  norm over
the annulus bounded by a constant (independent of $A$) times
$\eta$. Using a cutoff function $\rho$, the connection form can
be cut off to  define a connection $\rho A$ on $\m C^2$ which  is
trivial outside $B(0,4)$,  anti-self-dual (with respect to the
standard metric) on a neighbourhood of  $\b B(0,3)$, and has
self-dual curvature everywhere bounded in $C^0$ by a term of
order $\eta$;  bounds of the same order hold if the flat metric
is replaced by the standard conformally flat metric on $S^4$.
The following result shows that $\rho A$ can be perturbed into an
anti-self-dual connection on $S^4$:
\bigskip
\noindent {\bf
Lemma~2.3.} \quad  \sl Let $Y$ be a compact Riemannian
$4$-manifold with anti-self-dual Weyl curvature and positive
scalar curvature and let $E$ be a unitary bundle over $Y$.  Then
there  are constants $\delta, C>0$ such that if $A$ is a smooth
connection on  $E$ satisfying $\sup|F_+(A)| \le \delta $ it
follows that there is a solution $a \in \Lambda^1({\rm End}_0 E)$
to the equation $F_+(A+a)=0$ satisfying
$$
\norm{a}_{L^4} +
\norm{d_0 a}_{L^2} \le C\sup|F_+(A)|\;. \eqno(2.4)
$$
($d_0 a$ is
the full covariant derivative of $a$ using the connection $A$.)
\rm
\medskip
\noindent {\bf Proof: } The argument is modeled on
those in [FU, Ch.~7], using the Continuity Method to solve the
differential equation.  The solution $a$ will be of the form $a=
d_+^*b$ for $b\in \Lambda_+^2({\rm End}_0E)$, where the subscript
$A$ has been dropped to simplify the notation---throughout the
proof, all differential operators acting on bundle-valued forms
are those which are induced by the connection $A$ and the
Riemannian connection induced by the metric, with $d_0$  denoting
full covariant differentiation.

For $t \in [0,1]$  the equation $F_+(A+ a_t)=(1-t)F_+(A)$ can be
rewritten
$$
d_+^{}d_+^*b_t+[a_t,d_+^*b_t]_+ = -tF_+(A)\,, \quad
a_t=d_+^*b_t \;.   \eqno (2.5)
$$
The Weitzenb\"ock formula
(6.26) of [FU] (with the orientation reversed) reads
$$
d_+^{}d_+^* = d_0^*d_0^{} + R/3 + [F_+(A),\,\cdot\,] \eqno (2.6)
$$
where $R$ is the scalar curvature.  If $R_0 = \inf R$ and $F_+
:= F_+(A)$ then for $b \in \Lambda^2_+$
$$
\eqalign{ \langle
b,d_+^{}d_+^*b +[a_t,b]_+ \rangle &\ge  \langle d_0 b,d_0
b\rangle + {1\over 3}\langle b,Rb\rangle -
(\norm{a_t}_{L^2}+\norm{F_+}_{L^2})\norm{b}_{L^4}^2 \cr
&\ge
(c_1-2\norm{a_t}_{L^2}-\norm{F_+}_{L^2})\norm{b}_{L^4}^2 +
{R_0\over 6 }\norm{b}_{L^2}^2}
$$
where $c_1$ is the Sobolev
constant such that $c_1\norm{f}_{L^4}^2 \le \norm{df}_{L^2}^2
+R_0/6\norm{f}_{L^2}^2$ for $f\in L^2_1(Y)$.  Thus if
$\norm{a_t}_{L^2}+\norm{F_+}_{L^2} \le c_1/2$, the linearisation
of the operator $b \mapsto d_+^{}d_+^*b +[d_+^*b,d_+^*b]_+$ at
$b=b_t$ has no kernel and is therefore an isomorphism.

Next, using (2.6) as above yields
$$
\eqalign {{1\over
2}d^*d|b_t|^2 &= -|d_0 b_t|^2 + \langle b_t,d_0^*d_0^{}
b_t\rangle \cr
&\le -|d_0 b_t|^2 -{R_0\over 3}|b_t|^2 +
t|b_t||F_+| + t^2|b_t|^2|F_+| + |b_t||a_t|^2 \cr
&\le
-(1-|b_t|)|d_0 b_t|^2 - ({R_0\over 3}-\delta)|b_t|^2 +
|b_t||F_+|\;. }
$$
If $\sup |b_t| \le 1$ and $\delta \le R_0/6$
then at the maximum of $|b_t|$ the left-hand side of this
equation is non-negative giving the a priori bound
$$
\sup|b_t|
\le 6\sup|F_+|/R_0\;. \eqno(2.7)
$$
Moreover, if $|b_t| \le 1/2$,
the inequality above implies
$$
d^*d|b_t|^2 + |d_0 b_t|^2 +
{R_0\over 3}|b_t|^2 \le 2|b_t||F_+| \;. \eqno (2.8)
$$
Integrating both sides  also yields the a priori bound
$\norm{a_t}_{L^2}^2 \le \norm{d_0 b_t}^2_{L^2} \le
2\sup|b_t|\norm{F_+}_{L^1}$, so the invertiblity of the
linearized operator in (2.5) is guaranteed by a suitable bound on
$\sup|F_+|$ and the side condition $\sup|b_t| \le 1/2$.

\medskip
The Weitzenb\"ock formula (6.25) of [FU] for $1$-forms
is
$$
2d_+^*d_+^{} + dd^*=d_0^*d_0^{} + Ric -2[\,\cdot\,\lcorner F_-]
$$
where $Ric$ is the Ricci curvature of the metric.  Since $d^*a_t
=0$ and $d_+a_t = -tF_+ -a_t\wedge a_t$, it follows that
$$
{1\over 2}d^*d|a_t|^2 +|d_0 a_t|^2 =   -\langle
a_t,Ric(a_t)\rangle +2\langle a_t,[a_t\lcorner F_-]\rangle  -
2\langle a_t,d_+^*(tF_+ + (a_t\wedge a_t)_+) \rangle \;. \eqno
(2.9)
$$
Integrating both sides and estimating the right gives
$$
\eqalignno {\norm{d_0 a_t}_{L^2}^2 &\le
\norm{Ric}_{L^2}\norm{a_t}_{L^4}^2 + 2
\norm{a_t}^2_{L^4}\norm{F}_{L^2} + 4\norm{a_t}_{L^4}^2 +
4\norm{F_+}^2_{L^2} &\cr
&\le c_2\norm{a_t}_{L^4}^2
+4\norm{F_+}^2_{L^2}\;,  &(2.10) }
$$
where $c_2$ is a constant
depending only on the Riemannian metric (and the Chern classes of
$E$ assuming $\delta \le 1$.). To deal with the $L^4$ norm of
$a_t$  multiply (2.8) through by $|a_t|^2 $ and integrate. This
gives
$$
\eqalignno{ \int |a_t|^2|d_0 b_t|^2\,dV &\le 2\int |b_t|
|F_+||a_t|^2 \,dV - \int |b_t|^2 d^*d|a_t|^2\,dV &\cr
&= 2\int
|b_t| |F_+| |a_t|^2 \,dV +2\int |b_t|^2 |d_0 a_t|^2 \,dV - 2\int
|b_t|^2 K \,dV & (2.11) }
$$
where $K$ is the right-hand side of
(2.9). That is,
$$
-\int |b_t|^2 K\, dV =  \int
|b_t|^2\big[\langle a_t,Ric(a_t)\rangle -2\langle
a_t,[a_t\lcorner F_-]\rangle\big]\,dV   +\int|b_t|^2\big[\langle
a_t,d_+^*(tF_++(a_t\wedge a_t)_+) \rangle\big]\,dV \;.
$$
The
first integral on the right-hand side is bounded by
$c_2\sup|b_t|^2\norm{a_t}^2_{L^4}$;  the second is equal to
$$
\eqalign{ \hskip-.2in\int \langle d_+(|b_t|^2a_t),tF_+
+(a_t\wedge a_t)_+\rangle \,dV &\,=\,-\norm{|b_t|[tF_++(a_t\wedge
a_t)_+]}^2_{L^2}\cr
& \qquad +\int \langle d|b_t|^2\wedge
a_t,tF_++(a_t\wedge a_t)_+\rangle\,dV \cr
& \le 2\sup
|b_t|\norm{|d_0b_t||a_t|}_{L^2}^2(1+\norm{F_+}_{L^2}^2) }
$$
using $|a_t| \le |d_0 b_t|$ and Young's Inequality. Substituting
these bounds  back into (2.11) and applying (2.10) then gives
$$
\eqalign{ \int |a_t|^2|d_0 b_t|^2\,dV &\le
2\sup|b_t|\big[\norm{F_+}_{L^2}\norm{a_t}_{L^4}^2
+4\sup|b_t|(c_2\norm{a_t}^2_{L^4}+\norm{F_+}_{L^2}^2) \cr
&\qquad
\qquad + 2\norm{|d_0
b_t||a_t|}_{L^2}^2(1+\norm{F_+}_{L^2}^2)\big] \cr
& \le
4\sup|b_t|\big[c_3\norm{|d_0
b_t||a_t|}_{L^2}^2)+\norm{F_+}_{L^2}^2\big]   }
$$
where $c_3$
depends only on the geometry.  Rearranging terms gives
$$
(1-4c_3\sup|b_t|)\int|a_t|^2|d_0 b_t|^2\,dV \le
4\sup|b_t|^2\norm{F_+}^2_{L^2}
$$
so if $\sup|b_t| \le 1/8c_3$
then $\norm{a_t}_{L^4}^4 \le c_4\sup|F_+|^4$, using the bound on
$\sup|b_t|$ previously obtained.  Feeding this estimate back into
(2.10) then gives an a priori bound on $\norm{d_0 a_t}_{L^2}$ of
the form $c_5\sup |F_+|$, and the remainder of the proof is a
straight-forward application of the continuity method as in the
proof of Theorem~7.27 of [FU]. The details will be omitted. \qed
\bigskip

It follows from Lemma~2.3 that there is a sequence of
perturbations $\{a_i\}$ with $\norm{a_i}_{L^4}+\norm{d_{\rho
A_i,0}a_i}_{L^2}$ uniformly bounded by a constant of order $\sup
|F_+(\rho A_i)|$ such that $\rho A_i+a_i$ is an anti-self-dual
connection on an $SU(r)$ bundle over $S^4$, and the bound on
$\norm{a_i}_{L^2_1}$ implies that the transferred sequence is
bubbling the same amount of charge at $0\in \m C^2 \subset S^4$
as the original sequence was bubbling at $x_0$, $k>0$ units say.

\bigskip
\bigskip
\bigskip

\noindent{\bf 3. \quad Curvature of instantons on $S^4$.}

\bigskip
The Atiyah-Ward correspondence between instantons on
$S^4$ and monads on $\m{CP}_3 $ provides a description of the
perturbed sequence of transferred connections in terms of unitary
monads; that is, linear algebraic data.  To construct ``local"
monads for the {\it unperturbed\/} sequence requires good
estimates on the curvature of the instantons in terms of the
corresponding monads, and obtaining such estimates is the object
of this section.
\medskip
Recall the description of
anti-self-dual connections on $S^4$  in terms of holomorphic
bundles and monads on $\m P_2$ given in  [D1]: every holomorphic
$r$-bundle $E$ with $c_2=k$ which is trivial  on the line
$L_{\infty}$ at infinity is isomorphic to the cohomology of a
monad of the form
$$
\m{P}_2:\qquad 0 \longrightarrow K(-1)
\buildrel{A}\over\longrightarrow W
\buildrel{B}\over\longrightarrow K(1) \longrightarrow 0\quad
,\eqno(3.1)
$$
where $K$ and $W$  are hermitian vector spaces of
dimension  $k$ and $2k+r$ respectively and $A,\,B$ are linear
maps depending linearly on the  homogeneous coordinates $Z\in \m
P_2$: $A(Z)= A_0Z^0+A_1Z^1+A_2Z^2$  and similarly for $B$.  If
the coordinates are chosen so that  $L_{\infty}=
\{Z=(Z^0,Z^1,Z^2)\mid Z^2=0\}$ then bases for $K,\,W$ can  be
chosen so that if $R$ is the $r$-dimensional vector space
$R=(Im\,A_0+Im\,A_1)^{\perp}$ then  $W= K\oplus K \oplus R$,
$A(Z^0,Z^1,Z^2)=\bmatrix{Z^0{\bf 1}_K+ Z^2a_0\cr
Z^1{\bf
1}_K+Z^2a_1\cr
\hfill Z^2 a_2 } $ and $ B(Z^0,Z^1,Z^2)=[Z^1{\bf
1}_K+Z^2a_1,-Z^0{\bf 1}_K-Z^2a_0,Z^2b_2]$ for some  endomorphisms
$a_0,\,a_1$ of $K$ and some $a_2 \in  Hom(K,R),\; b_2 \in
Hom(R,K)$. The monad condition becomes  $ a_1a_0-a_0a_1 +
b_2a_2=0$ and the non-singularity  requirement is that $A(Z)$ be
injective and $B(Z)$ be surjective at  each point $Z\in \m P_2$.
The main result of [D1], proved using  geometric invariant
theory, is that each such monad is isomorphic  to  a
corresponding monad on  $\m P_3$  restricted to the plane
$Z^3=0$, where the latter monad possesses a unitary structure and
corresponds to a self-dual Yang-Mills connection (instanton) on
$S^4$.  This is an alternate formulation of the result that every
holomorphic bundle on $\m P_2$ which is trivial on $L_{\infty}$
admits a hermitian connection which is anti-self-dual with
respect  to the flat metric on $\m C^2$.  Expressed in terms of
the   monads above, the unitary structure is equivalent to the
condition that
$$
a_0^{}a_0^*-a_0^*a_0^{}+a_1^{}a_1^*-a_1^*a_1^{}
+b_2^{}b_2^*-a_2^*a_2^{}=0\;. \eqno(3.2)
$$

\medskip
Consider now the curvature of the induced unitary
connection on the bundle $E$ which is the cohomology of the monad
(3.1) (which is assumed to satisfy (3.2)).  Let $(z^0, z^1)$ be
inhomogeneous coordinates on $\m C^2=\{Z\in\m P_2\mid
Z^2\not=0\}$. Orthogonal projection $W\to E= \ker B \cap \ker
A^*$ is then given by $\pi_E= 1-A(A^*A)^{-1}A^*-B^*(BB^*)^{-1}B$,
from which  it follows that the curvature of the induced
connection on $E$ is $F=\pi_E[dA(A^*A)^{-1}\wedge
dA^*+dB^*(BB^*)^{-1}\wedge dB^*]\mid_E$. The unitary condition on
the monad implies that  $ A^*A=BB^* =:\psi\;, $  so if $i_E$
denotes the inclusion of $E$ into $W$ then the curvature is given
explicitly by
$$
F=\pi_E\bmatrix{\psi^{-1}&0&0\cr0&\psi^{-1}&0\cr
0&0&0}\bmatrix{dz^0\wedge d\bar z^0-dz^1\wedge d\bar z^1 &
2dz^0\wedge d\bar z^1 &0 \cr
2dz^1\wedge d\bar z^0 & dz^1\wedge
d\bar z^1-dz^0\wedge d\bar z^0 &0\cr
0&0&0}i_E \;.\eqno (3.3)
$$

If $p_{ij}$ is the $i,j$-th block in $\pi_E$, then (3.3) gives
$$
*\tr F^2
=8\tr[(p_{00}\psi^{-1}+p_{11}\psi^{-1})^2+
2(p_{00}\psi^{-1}p_{11}\psi^{-1}-
p_{01}\psi^{-1}p_{10}\psi^{-1})]  \le
16\tr[(p_{00}\psi^{-1}+p_{11}\psi^{-1})^2\,]\;;
$$
(the
combinatorial factor $8$ comes from the volume form
$(1/4)dz^0\wedge dz^1\wedge d\bar z^0\wedge d\bar z^1$).  Since
$\pi_E$ and $\psi$ are positive semi-definite with the former a
projection of rank $r$, $*\tr F^2 \le
16[\tr(p_{00}+p_{11})]^2[\tr\psi^{-1}]^2 = 16[r-\tr
p_{22}]^2[\tr\psi^{-1}]^2 \le 16r^2[\tr\psi^{-1}]^2$.  Thus
$$
|F| \le
4\tr(a_2^{}\psi^{-1}a_2^*+b_2^*\psi^{-1}b_2^{})\tr\psi^{-1}\le
4r\,\tr \psi^{-1}\;. \eqno(3.4)
$$

As the monad degenerates (that is, $a_0,a_1,a_2,b_2 \to 0$) the
automorphism $\psi\in End(K)$ approaches $|z|^2{\bf 1}_K$.  The
following lemma provides the estimates required for the next
section. In this lemma, the underlying Riemannian metric is the
standard flat one.
\bigskip
\noindent
\noindent {\bf Lemma~3.5.} \quad
{\sl Let $f := [\tr\psi^{-1}]^{-1/2}$ and $g := \log\det\psi$.
Then

\medskip
\itemitem{(a)} \quad $|df|^2 \le k$.
\itemitem{(b)}
\quad $\lap f^{-2} = 4\tr[\psi^{-1}(p_{00}+p_{11})\psi^{-1}] \ge
0$ \quad {\rm (where $p_{ij}$ is the $(i,j)$-th block of
$\pi_E$)}.
\itemitem{(c)} \quad $|dg|^2 \le 4kf^{-2}$.
\itemitem{(d)} \quad $\lap g \le -4f^{-2}$.
\itemitem{(e)} \quad
$f|d^2f|+f^2|d^2g|$ is uniformly bounded by a combinatorial
constant depending only on $k$. }

\medskip
\noindent {\bf
Proof: } (a) \quad Set $\alpha_i := z^i+a_i$ for $i=0,1$, so
$A(z)=[\alpha_0\;\alpha_1\;a_2]^{\rm T}$ and  $d\psi=
d(A^*\!A)=dA^*A+A^*dA=\sum_i(\alpha_i d\bar z^i+\alpha_i^*dz^i)$.
Since $d((\tr\psi^{-1})^{-1/2})=
(1/2)(\tr\psi^{-1})^{-3/2}\tr[\psi^{-1}(dA^*A+A^*dA)\psi^{-1}]$
it follows  $|d((\tr\psi^{-1})^{-1/2})|^2 =
(\tr\psi^{-1})^{-3}\sum_i |\tr(\psi^{-1}\alpha_i\psi^{-1})|^2 \le
(\tr(\psi^{-1}))^{-3}\tr(\psi^{-3})k \le k$.

\smallskip
\noindent (b) \quad $-d\cdot d \tr\psi^{-1}=d\cdot(\tr
\psi^{-1}(dA^*A+A^*dA)\psi^{-1}) = 8\tr\psi^{-2}
-2\tr[\psi^{-1}(dA^*A+A^*dA)\psi^{-1}\cdot
(dA^*A+A^*dA)\psi^{-1}]
=8\tr\psi^{-2}-4\sum_i\tr[\psi^{-1}(\alpha_i^*\psi^{-1}\alpha_i
+\alpha_i\psi^{-1}\alpha_i^*)\psi^{-1}] =  4\tr[\psi^{-1}({\bf
1}_K-\alpha_0^{}\psi^{-1}\alpha_0^*-
\alpha_1^*\psi^{-1}\alpha_1^{})\psi^{-1}
+\psi^{-1}({\bf
1}_K-\alpha_1^{}\psi^{-1}\alpha_1^*-
\alpha_0^*\psi^{-1}\alpha_0^{})\psi^{-1})]
=4\tr[\psi^{-1}(p_{00}+p_{11})\psi^{-1}]$.

\smallskip
\noindent (c) \quad $d\log\det \psi =
\tr[\psi^{-1}(dA^*A+A^*dA)]$, so  $|d\log\det \psi|^2=
\tr[\psi^{-1}(dA^*A+A^*dA)]\cdot \tr[\psi^{-1}(dA^*A+A^*dA)]=
4\sum_i|\tr\psi^{-1}\alpha_i|^2 \le 4k\tr\psi^{-1}$.

\smallskip
\noindent (d) \quad $-d\cdot d \log\det\psi =
-8\tr\psi^{-1} +\tr[\psi^{-1}(dA^*A+A^*dA)\psi^{-1}\cdot
(dA^*A+A^*dA)] =  -8\tr\psi^{-1} + \hfil\break
4\sum_i\tr[\psi^{-1}\alpha_i\psi^{-1}\alpha_i^*] =
-4\tr\psi^{-1} -2\tr[\psi^{-1}(p_{00} + p_{11})] \le
-4\tr\psi^{-1}$.

\smallskip
\noindent (e) \quad Straightforward calculation as
above. \qed

\bigskip
\noindent {\bf Lemma~3.6.} \quad {If $0\le p < 2$,
$q\ge 0$, and $c$  is a constant with $c\ge g+{2k(2q-2+p)\over
2-p}$ in a ball $B$ then
$$
\lap[f^{-p}(c-g)^q] \ge
2qf^{-p-2}(c-g)^{q-1} \quad \hbox{ in $B$.}
$$
}

\medskip
\noindent {\bf Proof:} \quad By direct calculation and  using
Lemma~3.5,
$$
\eqalign{ \lap [f^{-p}(c-g)^q] & =
\lap(f^{-p})(c-g)^q + f^{-p}\lap(c-g)^q  -2d(f^{-p})\cdot
d(c-g)^q\;, \cr
\lap(f^{-p}) &=
(p/2)f^{2-p}\lap(f^{-2})+p(2-p)f^{-p-2}|df|^2\cr
& \ge
p(2-p)f^{-p-2}|df|^2\;,\cr
\lap(c-g)^q & =-q(c-g)^{q-1}\lap g
-q(q-1)(c-g)^{q-2}|dg|^2 \cr
&\ge 4q(c-g)^{q-1}f^{-2}
-q(q-1)(c-g)^{q-2}|dg|^2\;. }
$$
By Young's Inequality,
$$
\eqalign{ 2d(f^{-p})\cdot d(c-g)^q &=2pqf^{-p-1}(c-g)^{q-1}df
\cdot dg \cr
&\le p(2-p)(c-g)^qf^{-p-2}|df|^2 + {pq^2\over
2-p}(c-g)^{q-2}f^{-p}|dg|^2}\,.
$$
Combining the above and using
Lemma~3.5(c) to estimate $|dg|^2$ gives $\lap[f^{-p}(c-g)^q]  \ge
4qf^{-p-2}(c-g)^{q-2}\big(c-g-{k(2q-2+p)\over 2-p}\big)$ from
which the desired inequality follows. \qed

\bigskip

\noindent {\bf Remark:} \quad Integrating the inequality of the
lemma over the unit ball shows that the function
$f^{-2}(c-g)^q=\tr \psi^{-1}(c-\log\det\psi)^q$ is uniformly
bounded in $L^p$ for any $p<2$ and $q \ge 0$ as
$|a_0|+|a_1|+|a_2|\to 0$.  This fails for general
$a_0,\,a_1,\,a_2$ not satisfying the reality constraint (3.2)
required to correspond to an instanton.
\bigskip

Suppose now that the monad (3.1) corresponds to an anti-self-dual
connection on a bundle over $S^4$ which has been obtained from
the procedure of  \S 2: a concentrated Hermitian-Einstein
connection $A$ on $(X,\omega)$ has been cut off and transferred
to $\rho A$ on $S^4$, then perturbed to an anti-self-dual
connection $\rho A +a$ on an $SU(r)$ bundle over $S^4$ with
$c_2=k$, with $a=d_{\rho A}^{+*}b$.   The cutoff function $\rho$
is identically $1$ on a neighbourhood of $\b B(0,3)$ and vanishes
outside $B(0,4)$.

The following lemma provides useful estimates on the size of the
perturbation $a$  in terms of the monad (3.1) and the
automorphism $\psi$:
\bigskip
\noindent {\bf Lemma~3.7.} \quad
\sl   There is a  constant $C_0$ independent of $A$ such that if
$\eta = \sup|F_+(\rho A)|$ is sufficiently small, then
$$
|a|^2
\le C_0^2\eta^2\tr\psi^{-1} \qquad\hbox{in $B(0,3)$.} \eqno(3.8)
$$
\medskip
\rm
\noindent {\bf Proof:} \quad  The underlying
Riemannian metric is the standard metric on $S^4$, and $d_0$ is
the full covariant derivative induced by $\rho A$ and this
metric.  The calculations of \S2 yield  bounds similar to (2.8)
of the form
$$
\lap |b|^2 + |d_0 b|^2 \le 2|b||F_+(\rho A)|,
\quad\qquad |b| \le c_1\eta \eqno(3.9)
$$
for some constant
$c_1$ independent of $A$.

\smallskip
The Weitzenb\"ock Formula (6.25) of [FU] gives
$\lap_0a=-Ric(a)-2d_+^*[F_+(\rho A)-(a\wedge a)_+] +2[a\lcorner
F_-(\rho A)]$, so
$$
\eqalignno{ \lap|a|^2 + 2|d_0a|^2 =2\langle
a,\lap_0a\rangle
&\le 4|a||d_+^*F_+(\rho A)|+4|F(\rho A)||a|^2 +8|a|^2|d_0a| &\cr
&\le 4|a||d_+^*F_+(\rho A)| +4|F(\rho A+a)||a|^2+
12|a|^2|d_0a|+4|a|^4 \;.&(3.10) }
$$

Fix a smooth cutoff function $\chi$  which is identically $1$ on
$B(0,3)$ and which is supported $\{z\in B(0,4)\mid \rho(z)=1\}$.
With $f := [\tr\psi^{-1}]^{-1/2}$ as in Lemma~3.5 it follows from
that lemma that both $|\lap(\chi^2f^2)|$, $|d(\chi f)|$ are
uniformly bounded by a fixed constant, so
$\lap[\chi^2f^2|a|^2]=\lap(\chi^2f^2)|a|^2+
\chi^2f^2\lap|a|^2-2\langle
d(\chi^2f^2),d|a|^2\rangle \le Const.(|a|^2+\chi
f|a||d_0a|)+\chi^2f^2\lap|a|^2$.  Since $F_+(\rho A) \equiv 0$ in
$B(0,3)$ and  $|F(\rho A+a)| \le Const.f^{-2}$ there,  (3.10)
implies that $ \lap[\chi^2f^2|a|^2] + 2\chi^2f^2|d_0a|^2 \le
Const.(|a|^2+\chi
f|a||d_0a|)+12|a|^2\chi^2f^2|d_0a|+4\chi^2f^2|a|^4 $. Applying
Young's inequality, it follows that there is a constant $c_2$
independent of $A$ such that
$$
\lap[\chi^2f^2|a|^2] +
\chi^2f^2|d_0a|^2 \le c_2(|a|^2+\chi^2f^2|a|^4)\; \eqno (3.11)
$$

If $Q := \sup(\chi f|a|)$ then $\lap[\chi^2f^2|a|^2] \le
c_2(1+Q)|a|^2$, so by (3.9) it follows that
$\lap[\chi^2f^2|a|^2+c_2(1+Q)|b|^2] \le 0$.  By the Maximum
Principle and (3.9), $\chi^2f^2|a|^2+c_2(1+Q)|b|^2 \le
c_2(1+Q)c_1^2\eta^2$ implying in particular that $Q^2\le
c_1^2c_2\eta^2(1+Q)$ and yielding the desired result. \quad \qed

\bigskip
\noindent {\bf Remark:}\quad By differentiating again
and expressing $d_0F(\rho A)$ in terms of $d_{\rho A+a,0}F(\rho A
+ a)$ (which is of order $f^{-3}$) and $d_0a$, it is
straightforward to show that the same methods yield a bound on
$|d_0a|$ of the form $|d_0a| \le Const.\eta \tr\psi^{-1}$ in
$B(0,2)$, and therefore there is a bound of the form $|F(\rho A)|
\le Const.\eta.\tr\psi^{-1}$ in $B(0,2)$.

\bigskip
\bigskip
\bigskip

\noindent{\bf 4. \quad Perturbation of monads I.}

\bigskip

Let $(A,\rho,a)$ be as at the end of \S3, and let $E$ be the
unitary $r$-bundle over $S^4$ on which $\rho A+a$ is defined. The
construction of a monad on $\m P_2$ (restricted from $\m P_3$)
corresponding to the instanton $\rho A + a$ can be viewed as a
particular embedding of the bundle $E$ into a larger
topologically trivial bundle which is equipped with a flat
connection so that the connection on $E$ is that which is induced
by Hermitian projection.  That is, splitting the bundle $W$ of
(3.1) using the Hermitian metric gives an isomorphism of $W$ with
$K(-1)\oplus E\oplus K(1)$ so that the connection $d_W$ on $W$ is
identified as
$$
\bmatrix{K(-1) \cr
E\cr K(1)} \owns \bmatrix{p\cr q\cr r} \;\mapsto\; \bmatrix{d_{K(-1)}p\cr d_Eq
\cr
d_{K(1)}r} +  \bmatrix{0 & \alpha & \beta\cr
-\alpha^* &0&\gamma \cr
-\beta^* &-\gamma^*&0} \bmatrix{p\cr q\cr r}  \eqno (4.1)
$$
where $\alpha \in \Lambda^{(0,1)}\otimes Hom(E,K(-1))$, $\beta
\in \Lambda^{(0,1)}\otimes Hom(K(1),K(-1))$ and $\gamma \in
\Lambda^{(0,1)}\otimes Hom(K(1),E)$ are the $\db$-closed
$(0,1)$-forms representing the extensions; i.e., the second
fundamental forms.

In this and the next section it will be shown that the
perturbation $a$ of the connection $\rho A$ is induced by a
corresponding small perturbation of the flat connection on $W$.
Since no curvature can bubble from a sequence of nearly-flat
connections on $W$, the degeneration of the connections on $E$ in
a neighbourhood of $0$ will be manifested as the degeneration of
the holomorphic maps $K(-1) \to W$, $W\to K(1)$, where $W$ is now
equipped with its new holomorphic structure.  This degeneration
is then amenable to analysis by algebraic methods.

Because the procedure for constructing the perturbation of the
monad is quite lengthy and somewhat indirect, the method is
outlined below.
\medskip

Let $d_W'$ be the connection on $W$ which is constructed from the
connection $\rho A$ on $E$ together with the connections on
$K(\pm 1)$ induced by $(4.1)$ using the same second fundamental
forms;  that is (abusing notation slightly), $d_W'=d_W-a$.  Thus
$$
F(d_W')=F(d_W)-d_W a+a\wedge a = \bmatrix{0 & -\alpha\wedge a
& 0 \cr
a\wedge \alpha^* & -d_{\rho A}a-a\wedge a & -a\wedge
\gamma \cr
0 & \gamma^*\wedge a &0 }\;. \eqno (4.2)
$$

The first perturbation of $d_W'$ which is sought will be obtained
from a perturbation
$(\alpha+\delta\alpha,\beta+\delta\beta,\gamma+\delta\gamma)$ of
$(\alpha,\beta,\gamma)$ so that the only non-zero $(0,2)+(2,0)$
component of the curvature is that coming from $F(\rho A)$. Using
(4.2), this means that with the connections on $K(\pm1)$ fixed
and with the connection $\rho A$ on $E$, the equations to be
satisfied are
$$
\eqalignno{ \db \delta\alpha  &=\alpha\wedge a''
&(4.3)(a)\cr
\db\delta\beta &= -[\alpha\wedge\delta\gamma +
\delta\alpha \wedge \gamma +\delta\alpha\wedge \delta\gamma] &(b)
\cr
\db\delta\gamma&=-a''\wedge\gamma &(c)\cr}
$$
where $a''$ is
the $(0,1)$ component of $a$.

The system (4.3) is non-linear only in an elementary way and can
be solved  via linear equations, first solving (a) and (c), and
then (b). The standard approach would be to take
$(\delta\alpha,\delta\beta,\delta\gamma)
=(\db^*x,\db^*y,\db^*z)$, where adjoints are with respect to the
Fubini-Study metric on $\m P_2$. For any of the $\db$-operators
involved in the system, the kernel of $\lap''$ on $\Lambda^{0,2}$
is canonically dual to the kernel of $\lap''$ on sections of the
dual bundle tensored with the canonical bundle ${\cal O}(-3)$ of
$\m P_2$.  Since $\rho A$ is flat outside  $B(0,4)$, any such
section of the latter on $\m P_2\backslash B(0,4)$ is a section
of ${\cal O}(-3)$ tensored with a trivial bundle twisted by
${\cal O}(1)$ or  ${\cal O}(2)$.  Such a section vanishes on all
lines not meeting $B(0,4)$ and hence everywhere, being the
solution of an elliptic equation.  Therefore, since $\lap''$ has
index $0$, there is no obstruction to solving (4.3)---indeed,
$(x,y,z)$ are uniquely determined.

Having thus perturbed the connection $d'_W$ to $d''_W$ say, the
second step is to remove the central part of the curvature:  an
inspection of the proof of the main result of [B4] shows that the
methods there remain valid in the current setting, even though
the connection is not everywhere integrable.  That is, there is a
solution $g\in Aut_{\m C}(W)$ on $\m C^2$ to the equation $\h
F(g\cdot d''_W)=0$ (where the $\Lambda$ operator
$\Lambda^{1,1}\to \Lambda^{0,0}$ is now that coming from the flat
metric on $\m C^2$), with $g^*g \to {\bf 1}$ along $L_{\infty}$.
The self-dual component of the curvature of this new connection
is now $gF^{0,2}(\rho A)g^{-1}+g^{*-1}F^{2,0}(\rho A)g^*$ (where
notation has again been abused by suppressing the inclusion map
of $E$ into $W$ and projection map $W\to E$).  To ensure that
this remains uniformly bounded and small therefore requires that
the automorphism  $g$ and its inverse remain uniformly bounded in
the annulus supporting $F_+(\rho A)$.  Such a bound can be
obtained from an appropriate bound on $\h F(d''_W)$, which in
turn follows from good estimates on the solutions of (4.3). The
problem outstanding is to provide these estimates.

On a compact K\"ahler surface, the Weitzenb\"ock formulae
relating the $\db$-Laplacian $\triangle''=\db\db^*+\db^*\db$
associated with a connection on a bundle $E$ with the full
covariant Laplacian $\lap_0$ are
$$
\lap_0 = \cases{2\triangle''+
i\h F \quad &on \quad $\Lambda^{0,0}(E)$ \cr
2\triangle''+2i\Lambda F_-\wedge -i\Lambda R\wedge \quad &on
\quad $\Lambda^{0,1}(E)$ \cr
2\triangle'' - i\h F -i\h R \quad
&on \quad$\Lambda^{0,2}(E)$}\eqno (4.4)
$$
where
$R\in\Lambda^{1,1}$ is the Ricci curvature and $F$ is the
curvature of the connection.  (The equations themselves are
easily deduced from the Hodge Identities, using the fact that on
$\Lambda^{0,q}$, the exterior derivative $\partial$ agrees with
the full covariant derivative $\partial_0$.)

The connections on $K(-1)$ and $Hom(K(1),-)$ induced from the
monad have curvatures $F$ such that $iF$ is negative, becoming
highly so in the ball $B(0,4)$ as the monad degenerates.  This
sign makes estimation of solutions to (4.3) less than
straightforward, and forces a modification to the strategy
outlined above.  The approach taken here is to solve (4.3) in
$B(0,3)$ using Dirichlet boundary conditions and the flat
Riemannian metric, for which the solutions are more easily
estimated. These solutions can then be cut off with small error
to yield a smooth almost-flat connection on $W$ over $S^4$ to
which the remainder of the method can be applied.
\medskip
With
$A$ and $B$ as in (3.1) the second fundamental forms in (4.1) can
be identified as
$$
\alpha =-(A^*\!A)^{-1/2}\db
A^*\!\mid_E\,,\quad  \beta =(A^*\!A)^{-1/2}A^*(\db
B^*)(BB^*)^{-1/2}\,, \quad \gamma=\pi_E(\db B^*)(BB^*)^{-1/2}
$$
where the adjoints are with respect to the Hermitian metric on
$W$ and a fixed metric on $K$, the line bundle ${\cal O}(-1)$
assumed to have been trivialised over the affine portion of $\m
P_2$.  Thus  $|\alpha|^2, |\beta|^2, |\gamma|^2 \le 4\tr
\psi^{-1}$, with $\psi=A(z)^*\!A(z)$ as before.

Consider first the equation (4.3)(a):  $\db\delta\alpha =
\alpha\wedge a''$, where the $\db$-operator is that which is
induced on $Hom(E,K(-1))$ by the connection $\rho A$ on $E$ and
that which is induced by Hermitian projection from the monad
(3.1) on $K(-1)$.  If $\db_u$ denotes the standard $\db$-operator
on the trivial bundle $K$ over $\m C^2$ (the subscript indicates
``untwisted"), then  $\db= \psi^{1/2}\circ \db_u \circ
\psi^{-1/2}$ so (4.3)(a) is equivalent to the equation
$$
\db_u(\psi^{-1/2}\delta\alpha)=\psi^{-1/2}\alpha\wedge
a''\;.\eqno (4.5)
$$
(The twisting by ${\cal O}(-1)$ has been
taken into account here:  $\db=\Psi^{1/2}\circ \db' \circ
\Psi^{-1/2}$ where $\Psi := A(Z)^*A(Z)/|Z|^2$ and $\db'$ is the
standard $\db$-operator on ${\cal O}(-1) $ over $\m P_2$.) As
indicated above,  the solution of this last equation is obtained
by solving the Dirichlet problem
$$
(\lap_u''\tau=)\quad
\db_u^{}\db_u^*\tau = \psi^{-1/2}\alpha\wedge a'' \qquad \hbox{
in $B(0,3)$, $\tau=0$ on $\d B(0,3)$} \eqno (4.6)
$$
for $\tau
\in \Lambda^{0,2}(Hom(E,K(-1))$.  Since all sections here are
smooth, standard linear elliptic PDE theory gives a smooth
solution of the equation, but the key information required are
estimates on such a solution. In the following lemma, $d_{u0}$
denotes the full covariant derivative induced from the flat
metric on $\m C^2$ together with $d_u$ on $K(-1)$ and $\rho A$ on
$E$.
\bigskip
\noindent {\bf Lemma~4.7.} \quad  \sl If
$|a_0|+|a_1|+|a_2| \le 1$ and $\eta$ is sufficiently small, the
solution $\tau$ of (4.6) satisfies
$$
\eqalignno{ |\tau|^2
&\;\le\;  C_1^2\eta^2 \tr\psi^{-1}&(4.8)\cr
\lap|\psi^{1/2}\tau|^2 + |\psi^{1/2}d_{u0}\tau|^2 &\;\le\;
C_1^2\eta^2\tr\psi^{-1} &(4.9) \cr
|\psi^{1/2}\tau|^2
&\;\le\;C_1^2\eta^2(3k-\log\det\psi) &(4.10) }
$$
where $C_1$ is
a constant independent of $A$.
\medskip
\rm

\noindent {\bf Proof:} \quad  As in \S3, let $f=
[\tr\psi^{-1}]^{-1/2}$ and $g = \log\det\psi$. The third
Weitzenb\"ock formula of (4.4) with $\h F=0=\h R$ gives
$$
\lap|\tau|^2+2|d_{u0}\tau|^2\le
2|\tau||\lap_{u0}\tau|=4|\tau||\psi^{-1/2}\alpha\wedge
a''|\;.\eqno(4.11)
$$
By Lemma~3.5, the bound on $|\alpha|$ and
the Cauchy-Schwarz inequality, it follows that $\lap|\tau| \le
|\lap_0\tau| \le 2|\psi^{-1/2}\alpha\wedge a''| \le 2c_1\eta
f^{-3}$ for some constant $c_1$ independent of the connection
$A$.

The assumption  $|a_0|+|a_1|+|a_2| \le 1$ implies that $g\le 3k$
in $B(0,3)$.  Taking $c=5k$ and $p=1=q$ in Lemma~3.6  gives
$\lap[f^{-1}(5k-g)]\ge 2f^{-3}$ in $B(0,3)$, so by the Maximum
Principle it follows that
$$
|\tau| \le c_1\eta
f^{-1}(5k-g)\;.\eqno (4.12)
$$

To sharpen this estimate,  observe  that
$$
|\psi^{-1/2}\alpha|^2=\tr [\pi_E dA\psi^{-2}\cdot
dA^*]=2\tr[\psi^{-2}(p_{00}+p_{11})] = (1/2)\lap\tr\psi^{-1}\;,
\eqno (4.13)
$$
(where $p_{ij}$ is the $i,j$-th block in $\pi_E$
with respect to the decomposition $W = K \oplus K\oplus \m C^r$).
If $\lambda>0$ is a constant (to be fixed later) it follows from
(4.11) that
$$
\eqalignno{ \lap|\tau|^2+2|d_{u0}\tau|^2 &\le
4|\tau||a''||\psi^{-1/2}\alpha|& \cr
&\le
2\lambda^{-1}\eta^2|\psi^{-1/2}\alpha|^2+
2\lambda\eta^{-2}|\tau|^2||a''|^2
&\cr
&=
\lambda^{-1}\eta^2\lap\tr\psi^{-1}+
\lambda\eta^{-2}|\tau|^2|a|^2\;.
&(4.14) }
$$

As in (3.9), there is an inequality of the form
$$
\lap|b|^2 +
|d_0b|^2 \le 2|b||F_+(\rho A)| +c_2|b|^2 =c_2|b|^2\quad \hbox{ in
$B(0,3)$ }
$$
for some uniform constant $c_2$ independent of $A$,
this additional term arising from the fact that it is now the
flat metric on $\m C^2$ which is being used rather than the
conformally flat metric on $S^4$; note that implies $|a|^2
=|*_{S^4}db|^2 \le (1+|z|^2)^4|d_0b|^2$.  Using Young's
inequality,
$$
\eqalign{ \lap(|b|^2|\tau|^2) &=
\lap(|b|^2)|\tau|^2+|b|^2\lap|\tau|^2-2d|b|^2\cdot d|\tau|^2 \cr
& \le -|d_0b|^2|\tau|^2-2|b|^2|d_{u0}\tau|^2
+|b|^2[c_2|\tau|^2+\lambda^{-1}\eta^2\lap\tr\psi^{-1}+
\lambda\eta^{-2}|\tau|^2|a|^2]
+8|b||\tau||d_0b||d_{u0}\tau| \cr
& \le
-(1/2)|d_0b|^2|\tau|^2+30|b|^2|d_{u0}\tau|^2
+|b|^2[c_2|\tau|^2+\lambda^{-1}\eta^2\lap\tr\psi^{-1}+
\lambda\eta^{-2}|\tau|^2|a|^2]\,.
}
$$
If $c_3$ is a fixed uniform constant such that $|b| \le
c_3\eta$ in $B(0,3)$ then  from (4.14),
$$
\lap[|b|^2|\tau|^2+16c_3^2\eta^2|\tau|^2] \le  -(1/2-32\lambda
c_3^2)|d_0b|^2|\tau|^2-2c_3^2\eta^2|d_{u0}\tau|^2
+c_2c_3^2\eta^2|\tau|^2+
17\lambda^{-1}c_3^2\eta^4\lap\tr\psi^{-1}\;,
$$
so if $\lambda=1/(128c_3^2)$ it follows that
$$
\lap[|b|^2|\tau|^2+16c_3^2\eta^2|\tau|^2]
+(1/4)|d_0b|^2|\tau|^2+2c_3^2\eta^2|d_{u0}\tau|^2 \le
c_2c_3^2\eta^2|\tau|^2+17\lambda^{-1}c_3^2
\eta^4\lap\tr\psi^{-1}\;.
\eqno(4.15)
$$
By (4.12) and Lemma~3.6, $|\tau|^2 \le
c_1^2\eta^2f^{-2}(5k-g)^2 \le \lap[c_1^2\eta^2(7k-g)^3/6]$,  so
by the Maximum Principle it follows that
$|b|^2|\tau|^2+16c_3^2\eta^2|\tau|^2 \le
17\lambda^{-1}c_3^2\eta^4\tr\psi^{-1}+c_1^2c_3^2\eta^4(7k-g)^3/6$.
The hypotheses imply that $(7k-g)^3 $ is bounded in  $B(0,3)$ by
a uniform constant times $\tr\psi^{-1}$, implying in particular
that $|\tau|^2 \le Const.\eta^2\tr\psi^{-1}$ for some constant
independent of $A$, proving (4.8).
\smallskip
To prove (4.9) and
(4.10), $\lap|\psi^{1/2}\tau|^2=\lap\langle \tau,\psi\tau\rangle
= 2Re\langle \tau,\psi\lap_{u0} \tau\rangle -4Re \langle
d_{u0}\tau, (d\psi)\tau\rangle
-8|\tau|^2-2|\psi^{1/2}d_{u0}\tau|^2$.  By the Weitzenb\"ock
formula, $\lap_{u0}\tau= 2\lap_u''\tau=2\psi^{-1/2}\alpha\wedge
a''$, so
$$
\lap|\psi^{1/2}\tau|^2 + |\psi^{1/2}d_{u0}\tau|^2
\le 4Re\langle \tau,\psi^{1/2}\alpha\wedge a''\rangle +8|\tau|^2
\le Const.\eta^2\tr\psi^{-1}\;,
$$
where the second inequality
follows from the last bound on $|\tau|$ and Lemma~3.8. The bound
(4.10) on $|\psi^{1/2}\tau|$  then follows  from the Maximum
Principle and Lemma~3.5(c). \quad \qed

\bigskip

Equation (4.3)(c) can be solved in $B(0,3)$ in precisely the same
way as above;  the solution $\delta\gamma$ takes the form
$\delta\gamma=(\db_u^*\mu)\psi^{1/2}$ for some $\mu \in
\Lambda^{0,2}(Hom(K(1),E)$ satisfying the same estimates as
$\tau$ in Lemma~4.7;  (more precisely, $\mu^*$ does).  Equation
(4.3)(b) is however more delicate:  using the untwisted
connection on $Hom(K(1),K(-1))$ it takes the form
$$
\db_u(\psi^{-1/2}\delta\beta\psi^{-1/2}) =
-\psi^{-1/2}[\alpha\wedge\delta\gamma + \delta\alpha \wedge
\gamma +\delta\alpha\wedge \delta\gamma]\psi^{-1/2}
$$
where
$\delta\alpha = \psi^{1/2}\db_u^*\mu$ and
$\delta\gamma=(\db_u^*\nu)\psi^{1/2}$. As above the equation can
be solved in $B(0,3)$ by taking
$\delta\beta=\psi^{1/2}\db_u^*\nu\psi^{1/2}$ for  $\nu \in
\Lambda^{(0,2)}(Hom(K(1),K(-1)))$ vanishing on $\d B(0,3)$, but
complications arise from the fact that the right-hand side of the
equation
$$
\lap''_u\nu= -[\psi^{-1/2}\alpha\wedge\db_u^*\mu +
\db_u^*\tau \wedge \gamma\psi^{-1/2} +\db_u^*\tau\wedge
\db_u^*\mu] \eqno(4.16)
$$
is of order $f^{-4}$ making
application of Lemma~3.6 impossible.  However, using (4.13) as
was done in the proof of Lemma~4.7 overcomes this problem.
\bigskip
\noindent {\bf Lemma~4.17.} \quad \sl If
$|a_0|+|a_1|+|a_2| \le 1$ and $\eta$ is sufficiently small, the
solution $\nu$ of (4.16) satisfies
$$
\eqalignno{ |\nu|^2
&\;\le\;  C_2^2\eta^2 f^{-2}&(4.18)\cr
\lap\big[\,|\psi^{1/2}\nu\psi^{1/2}|^2+c_0(5k-g)(|\psi^{1/2}\tau|^2
+|\mu\psi^{1/2}|^2)\,\big]+ |\psi^{1/2}d_{u0}\nu\psi^{1/2}|^2
&\;\le\; C_2^2\eta^2(5k-g)^2f^{-2}&(4.19) \cr
|\psi^{1/2}\nu\psi^{1/2}|^2 &\;\le\;C_2^2\eta^2(5k-g)^2 &(4.20) }
$$
where $c_0,C_2$ are constants independent of the connection
$A$, and $f^{-2} = \tr\psi^{-1}$ and $g=\log\det\psi$.
\medskip
\rm
\noindent {\bf Proof:} \quad  As in the proof of Lemma~4.7,
$$
\lap|\nu| \le |\lap_{u0}\nu| =2|\lap_{u}''\nu| \le \eta
(|\psi^{-1/2}\alpha|^2+|\gamma\psi^{-1/2}|^2)+
2\eta^{-1}(|\db_u^*\tau|^2+|\db_u^*\mu|^2)\;.
$$
{}From (4.15) and
the corresponding inequality for $\mu$ it follows that
$$
\lap[|\nu|+
\eta^{-3}c_3^{-2}(|b|^2+16c_3^2\eta^2)(|\tau|^2+|\mu|^2)] \le
c_4[\eta^{-1}(|\tau|^2+|\mu|^2)+\eta\lap\tr\psi^{-1}]
$$
for some
uniform constant $c_4$, so by (4.8), (4.13) and the Maximum
Principle it follows as above that $|\nu| \le
c_5\eta\tr\psi^{-1}$.
\smallskip

Inequality (4.19) will follow once (4.20) has been established:
$$
\eqalign{ \lap|\psi^{1/2}\nu\psi^{1/2}|^2 &= \lap\langle
\nu,\psi\nu\psi\rangle \cr
&=
\langle\lap_{u0}\nu,\psi\nu\psi\rangle + \langle
\nu,\psi(\lap_{u0}\nu)\psi\rangle -8\langle
\nu,\psi\nu+\nu\psi\rangle - 2\langle
d_{u0}\nu,\psi(d_{u0}\nu)\psi\rangle\cr
& \quad -2\langle
d_{u0}\nu,(d\psi)\nu\psi+\psi\nu d\psi\rangle -2\langle
(d\psi)\nu,(d_{u0}\nu)\psi\rangle -2\langle \nu d\psi,\psi
d_{u0}\nu\rangle -2 \langle (d\psi)\nu,\nu d\psi\rangle\,. }
$$
Setting $X := \psi^{1/2}d_{u0}\nu\psi^{1/2}$, $Y :=
\psi^{-1/2}(d\psi)\psi^{-1/2}$ and $Z :=
\psi^{1/2}\nu\psi^{1/2}$, the terms involving a factor of $-2$ in
this equation are $-2[|X|^2+\langle X,YZ+ZY\rangle + \langle YZ +
ZY,X\rangle + \langle YZ,ZY \rangle] \le 2|ZY|^2+2|ZY|^2 +
2\langle ZY,YZ \rangle]$.  Since $|ZY|^2 =
\tr\psi^{-1/2}d\psi\nu^*\psi\nu \cdot d\psi\psi^{-1/2}= 2\tr
[\nu^*\psi\nu({\bf 1}_K-p_{00}+{\bf 1}_K-p_{11})] \le 4\langle
\nu,\psi\nu\rangle$ and similarly $|YZ|^2 \le 4\langle
\nu,\nu\psi\rangle$, it follows that
$$
\eqalign{
\lap|\psi^{1/2}\nu\psi^{1/2}|^2 &\le
\langle\lap_{u0}\nu,\psi\nu\psi\rangle + \langle
\nu,\psi(\lap_{u0}\nu)\psi\rangle+ 2\langle \nu
d\psi,(d\psi)\nu\rangle \cr
&\le
4|\psi^{1/2}\nu\psi^{1/2}||\psi^{1/2}(\lap_u''\nu)\psi^{1/2}|
+2|\psi^{1/2}\nu\psi^{1/2}||\psi^{-1/2}d\psi||\psi^{-1/2}d\psi||\nu|
\cr
& \le
4|\psi^{1/2}\nu\psi^{1/2}|\,[\,|\alpha\wedge\delta\gamma +
\delta\alpha \wedge \gamma +\delta\alpha\wedge
\delta\gamma|+c_6\eta f^{-2}]\cr
&\le
4Q(5k-g)[\eta^{-1}|\delta\alpha|^2+\eta^{-1}|\delta\gamma|^2+c_7\eta
f^{-2}] \quad \hbox{ for $Q :=
\sup\,\displaystyle{|\psi^{1/2}\nu\psi^{1/2}|\over 5k-g}$,} }
$$
where $c_6,c_7$ are constants independent of $A$.  From (4.9),
(4.10) and Lemma~3.5,
$$
\lap\big[(5k-g)|\psi^{1/2}\tau|^2\big] +
{1\over 2}(5k-g)|\psi^{1/2}d_{u0}\tau|^2  \le
c_8\eta^2f^{-2}(5k-g)\;,
$$
so
$$
\lap[|\psi^{1/2}\nu\psi^{1/2}|^2+8Q\eta^{-1}(5k-g)(|\psi^{1/2}\tau|^2
+|\mu\psi^{1/2}|^2)] \le 16Qc_8\eta f^{-2}(5k-g)\;.
$$
By Lemma~3.6 again,
$$
\lap[|\psi^{1/2}\nu\psi^{1/2}|^2+8Q\eta^{-1}(5k-g)(|\psi^{1/2}\tau|^2
+|\mu\psi^{1/2}|^2) -4Qc_8\eta (5k-g)^2] \le 0\;,
$$
implying in
particular that $|\psi^{1/2}\nu\psi^{1/2}|^2 \le
4Qc_8\eta(5k-g)^2$. Hence $|\psi^{1/2}\nu\psi^{1/2}|^2(5k-g)^{-2}
\le 4Qc_8\eta$, giving $Q \le 4c_8\eta $; that is,
$$
|\psi^{1/2}\nu\psi^{1/2}| \le 4c_8\eta(5k-g)\;.
$$
Note that this
implies $|\nu\psi^{1/2}|,\, |\psi^{1/2}\nu|\le 4 c_8\eta
(5k-g)f^{-1}$.

Returning now to the proof of (4.19),
$$
\eqalign{
\lap|\psi^{1/2}\nu\psi^{1/2}|^2 +
|\psi^{1/2}d_{u0}\nu\psi^{1/2}|^2 &\le
4|\psi^{1/2}\nu\psi^{1/2}|\,|\alpha\wedge\delta\gamma +
\delta\alpha \wedge \gamma +\delta\alpha\wedge \delta\gamma| +
c_9(|\nu\psi^{1/2}|^2 + |\psi^{1/2}\nu|^2) \cr
& \le
c_{10}[\eta(5k-g)|\alpha\wedge\delta\gamma + \delta\alpha \wedge
\gamma +\delta\alpha\wedge \delta\gamma|+
\eta^2(5k-g)^2f^{-2}]\cr
&\le
c_{11}[(5k-g)(|\delta\alpha|^2+
|\delta\gamma|^2)+\eta^2(5k-g)^2f^{-2}]\;,
}
$$
and therefore, as above it follows that
$$
\lap[|\psi^{1/2}\nu\psi^{1/2}|^2+2c_{11}(5k-g)(|\psi^{1/2}\tau|^2
+|\mu\psi^{1/2}|^2)]+ |\psi^{1/2}d_{u0}\nu\psi^{1/2}|^2 \le
c_{12}\eta^2(5k-g)^2f^{-2}\,,
$$
as required.  \quad \qed
\bigskip
\medskip

\noindent{\bf Remark:} \quad Using the fact that the curvature of
$\rho A$ is bounded by a constant times $\tr\psi^{-1}$ as
remarked at the end of \S3, the same methods as given here yield
bounds on $|\delta \alpha|, |\delta\gamma|, |\delta\beta|$ over
$\bar B(0,2)$ of the form $Const.\eta f^{-1}(c-g)^q$ for some
(fixed) $q$;  indeed, the majority of the work has already been
done.  Furthermore, it is easily shown that $|d_0\delta \alpha|,
|d_0\delta\gamma|, |d_0\delta\beta|$ are uniformly bounded in
$L^2(B(0,2))$ by a term of order $\eta$.

\bigskip
\bigskip
\bigskip

\noindent{\bf 5. \quad Perturbation of monads II.}

\bigskip

Consider now the curvature of the perturbed connection $d''_W$ on
$W$ over the ball $B(0,3)$ constructed in the last section: {
\def\a{\alpha}
\def\b{\beta}
\def\g{\gamma}
\def\e{\delta}
\def\w{\wedge}
$$
\bmatrix{
(\!-\!\a\!\w\!\e\a^*\!-\!\b\!\w\!\e\b^*\!-\!\e\a\!\w\!\a^* &
(\!d\e\a \!-\!\a\!\w\!a\!\!-\!\b\!\w\!\e\g^* &
(d\e\b\!+\!\a\!\w\!\e\g \cr
\!-\!\e\b\!\w\!\b^*\!-\!\e\a\!\w\!\e\a^*\!-\!\e\b\!\w\!\e\b^*) &
\!-\!\e\b\!\w\!\g^*\!-\!\e\b\!\w\!\e\g^*) &
\!+\!\e\a\!\w\!\g\!+\!\e\a\!\w\!\e\g) \cr
&&\cr
(-\!d\e\a^*\!+\!a\!\w\!\a^* \!-\!\e\g\!\w\!\b^*&  (\!-\![d_{\rho
A}a\!+\!a\!\w\! a]\!-\!\a^*\!\w\!\e\a\!-\!\g\!\w\!\e\g^* &
(d\e\g\!-\!a\!\w\!\g\!-\!\a^*\!\w\!\e\b \cr
\!-\!\g\!\w\!\e\b^*\!-\!\e\g\!\w\!\e\b^*) &
\!-\!\e\a^*\!\w\!\a\!-\e\g\w\g^*\!-
\!\e\a^*\!\w\!\e\a\!-\!\e\g\!\w\!\e\g^*)
& \!-\!\e\a^*\w\b\!-\!\e\a^*\!\w\!\e\b) \cr
&&\cr
(\!-\!d\e\b^*\!+\!\e\g^*\!\w\!\a^*&
(d\e\g^*\!+\!\g^*\!\w\!a\!-\!\e\b^*\!\w\!\a&
(\!-\!\g^*\!\w\!\e\g\!-\!\e\g^*\!\w\!\g\!-\!\b^*\!\w\!\e\b \cr
\!+\!\g^*\!\w\!\e\a\!+\!\e\g^*\!\w\!\e\a^*) &
\!-\!\b^*\!\w\!\e\a\!-\!\e\b^*\!\w\!\e\a) &
\!-\!\e\b^*\!\w\!\b\!-\!\e\b^*\!\w\!\e\b\!-\!\e\g^*\!\w\!\e\g) }
$$
}
By construction, the $(0,2)$ and $(2,0)$ components of
$F(d''_W)$ are zero. Since both $\rho A$ and $\rho A +a$ are
anti-self-dual in $B(0,3)$ it follows that $|\h F(d''_W)| \le
C[\eta^{-1}(|\delta\alpha|^2+|\delta\beta|^2+|\delta\gamma|^2)+\eta
f^{-2}]+2(|\Lambda\d\delta\alpha|+|\Lambda\d\delta\beta|+
|\Lambda\d\delta\gamma|)$,
where the $\d$-operators here are the respective twisted
connections induced from the monad (3.1) and $\rho A$, and $C$ is
a combinatorial constant.

To estimate the terms involving derivatives of the perturbations,
observe that
$$
i\Lambda\d\delta\alpha=
i\Lambda\psi^{-1/2}\d_u(\psi^{1/2}\psi^{1/2}\db_u^*\tau)=
-\psi^{-1/2}\db_u^*(\psi\db_u^*\tau)
=i\Lambda(\psi^{-1/2}\d\psi\psi^{-1/2}\wedge\delta\alpha)\,,
$$
so  it follows that
$$
|\Lambda\d\delta\alpha|\le
|\psi^{-1/2}\d\psi\psi^{-1/2}||\delta\alpha|  \le
\eta|\psi^{-1/2}\d\psi\psi^{-1/2}|^2 + \eta^{-1}|\delta\alpha|^2
\le 2\eta\,\tr\psi^{-1}+ \eta^{-1}|\delta\alpha|^2\,.
$$
Similarly, $i\Lambda\d\delta\gamma=
-i\Lambda\delta\gamma\wedge\psi^{-1/2}\d\psi\psi^{-1/2}$ and
$i\Lambda \d\delta\beta=
i\Lambda[\psi^{-1/2}\d\psi\psi^{-1/2}\wedge\delta\beta-
\delta\beta\wedge\psi^{-1/2}\d\psi\psi^{-1/2}]$, so the same type
of estimate holds for  $|\Lambda\d\delta\beta|$ and
$|\Lambda\delta\gamma|$.  Thus for some new constant $C$,
$$
|\h
F(d''_W)| \le
C[\eta^{-1}(|\delta\alpha|^2+|\delta\beta|^2+|\delta\gamma|^2)+\eta
f^{-2}]\;.\eqno(5.1)
$$

Now fix a smooth cutoff function $\chi$ compactly supported in
$B(0,3)$ which is identically $1$ on a neighbourhood of $\bar
B(0,2)$, and replace $\tau, \mu$ and $\nu$ by $\chi\tau, \chi\mu$
and $\chi\nu$ respectively. Equations (4.9) and (4.19) remain
true if the constants $C_1$ and $C_2$ are altered, and if $\eta$
is sufficiently small (or the connection $\rho A$ is sufficiently
concentrated) the resulting perturbation to $d''_W$ is small in
$C^1$, bounded by a term of order $\eta\sup_{B(0,3)\backslash
B(0,2)}|\alpha|$.

Consider the connection on $W$ over $\m C^2$ obtained by cutting
off $\tau$, $\mu$ and $\nu$ in this way and replacing $\rho A$ by
$\rho A +(1-\rho)a$.  Outside  $B(0,4)$, it agrees with the
original flat connection $d_W$ on $W$; the $(0,2)$ and $(2,0)$
components of the curvature are supported in the annulus
$B(0,4)\backslash \bar B(0,2)$ and are bounded in $C^0$ by a term
of order $\eta$; and on a neighbourhood of $\bar B(0,2)$ it
agrees with $d''_W$ and therefore induces the connection  $\rho
A$ on $E$ by hermitian projection.  The notation $d''_W$ will be
retained for this new connection.

As noted in the introduction to the previous section, the methods
of [B4] can be applied to $d''_W$ to give a smooth solution $g\in
Aut_{\m C}(W)$ to the equation $\h F(g \cdot d''_W)=0$ with  with
$\det g\equiv 1$ and  $g^*g \to {\bf 1}_W$ on the line
$L_{\infty}$ at infinity in $\m P_2$.  Then  $\lap\log
|g|^2=2i\Lambda\dbd\log|g|^2 \le 2|\h F(d''_W)|$, so from (5.1),
(4.9) and (4.19) (with modified constants to account for the
introduction of $\chi$) it follows that in a neighbourhood of
$\bar B(0,3)$,
$$
\lap[\log|g|^2+2C\eta^{-1}(1+c_0(5k-\log\det\psi))
(|\psi^{1/2}\tau|^2+|\mu\psi^{1/2}|^2)+
|\psi^{1/2}\mu\psi^{1/2}|^2] \le C_3\eta
f^{-2}(5k-\log\det\psi)^2 \eqno(5.2)
$$
where $C_3$ is a constant
independent of $\rho A$.  Applying Lemma~3.6 with $p=0,\,q=3$
gives $\lap((7k-\log\det\psi)^3) \ge 6f^{-2}(7k-\log\det\psi)^2$
in $B(0,3)$.  If $(7k-g)_+$ denotes the positive part of  the
function $7k-g$ then for some suitable constant $C_4$ it follows
$$
\lap[(7k-\log\det\psi)_+^3+C_4(1+|z|^2)^{-1}] \ge
6f^{-2}(7k-\log\det\psi)_+^2 \qquad \hbox{ in $\m C^2$.}
$$
Therefore, by the Maximum Principle and the fact that $g^*g \to
1$ as $|z|\to\infty$ it  follows that
$$
\log|g|^2 \le \log
r+C_5\eta\,[(7k-\log\det\psi)_+^3+(1+|z|^2)^{-1}]
$$
for some
uniform constant $C_5$. Since $\det g=1$,  a corresponding
uniform bound exists for $g^{-1}$, and since $F_+(g\cdot
d''_W)=gF^{0,2}(d''_W)g^{-1}+{g^*}^{-1}F^{2,0}(d''_W){g^*}$ it
follows that $|F_+(g\cdot d''_W)| \le Const.\eta$.  Topological
triviality of $W$ implies $\norm{F(g\cdot d''_W)}_{L^2} \le
Const.\eta$ and therefore if $\eta$ is sufficiently small it
follows from Uhlenbeck's theorem [U1], [S] that after gauge
transformations, (a subsequence of) the connections $g\cdot
d''_W$  converges in $C^2$ (say) on $\m C^2$ to a smooth
connection on $W$ which is anti-self-dual in the complement of
the annulus $B(0,4)\backslash\bar B(0,2)$.

\smallskip
Since the sequence of connections  $g\cdot d''_W$ is
converging in $C^2$ and each defines an integrable connection in
a neighbourhood of $\bar B(0,2)$, there is a converging sequence
of holomorphic trivialisations of each of these connections in
this neighbourhood. Equivalently, there is a convergent sequence
of smooth automorphisms $g'$  of $W$ there such that
$g'\cdot(g\cdot d''_W)=d_W$. Using a cutoff function supported in
$B(0,3)$ which is the identity on $B(0,2)$, extend $g'$ to $\m
C^2$ and combine the product of automorphisms $g'g$ into one,
also denoted by $g$.  Thus $ g \cdot d''_W = d_W$ in $B(0,2)$,
$g^*g\to 1$ along $L_{\infty}$, $\det g=1$ (without loss of
generality) and the sequence of such automorphisms is converging
uniformly in $C^3$ (say) on compact subsets of $\m
C^2\backslash\{0\}$.

\smallskip
As the sequence of connections degenerates, the
cohomology of the monad (3.1) converges to the trivial rank $r$
bundle $R$ on compact subsets of $\m C^2\backslash \{0\}$, and
moreover the second fundamental forms $\alpha, \,\beta,\,\gamma$
converge to zero on such subsets---equations (4.6) and (4.16)
together with a priori $L^2$ estimates from (4.12) and (4.18)
imply that the perturbations
$\delta\alpha,\,\delta\beta,\,\delta\gamma$ also converge to $0$
uniformly on compact subsets of $B(0,3)\backslash\{0\}$.

The perturbation $a$ of the cut-off connection $\rho A$ satisfies
$[d_{\rho A + a}a-a\wedge a]_+ =0$ inside $B(0,3)$ and in view of
the  a priori $L^2_1$ bounds on $a$, it can be assumed to
converge weakly in $L^2_1(\m C^2)$  and uniformly in $C^1$ on
compact subsets of $\m C^2\backslash\{0\}$ to a form $a_{\infty}$
satisfying $[d_Ra_\infty -a_{\infty}\wedge a_{\infty}]_+=0$.
Regarding $a_{\infty}$ as defining a connection on $R$, the
Removable Singularities Theorem~[U2] gives a gauge in which
$a_{\infty}$ extends smoothly across the origin.  Combining the
Newlander-Nirenberg theorem with this gauge transformation yields
an  automorphism $g'_{\infty}$   of $R$ with unit determinant
which is smooth away from the origin such that
$a_{\infty}^{(0,1)} = (\db_Rg'_{\infty})g_{\infty}'^{-1}$.  Away
from $0$, extend $g'_{\infty}$ to all of $W$ by setting
$g'_{\infty}\mid_{R^{\perp}} = {\bf 1}_{R^{\perp}}$, so
$g'_{\infty}\cdot d''_{W,\infty} = d_W$ in $B(0,2)$.  If
$g_{\infty}$ is the limit of the sequence  of automorphisms
constructed in the previous paragraph, then $g_{\infty}\cdot
d''_{W,\infty} = d_W$ in $B(0,2)\backslash\{0\}$ also, so
$g_{\infty}^{}{g'}_{\infty}^{-1}$ is a holomorphic automorphism
of $W$ with unit determinant over $B(0,2)\backslash\{0\}$ which
extends smoothly  across the origin by Hartogs' Theorem.  Using a
cutoff function which is the identity on $B(0,3/2)$ and is
compactly supported in $B(0,2)$, extend $h_{\infty}$ smoothly to
all $\m C^2$ so that it has unit determinant everywhere and is
the identity outside $B(0,2)$, and then replace each of the
automorphisms $g$ of the sequence by $h_{\infty}^{-1}g$.  Then
inside $B(0,3/2)$, the identity $g\cdot d''_W =d_W$ is still
valid, each automorphism $g$ still has unit determinant, $g$ is
converging in $C^2$ on compact subsets of $\m
C^2\backslash\{0\}$, $g^*g \to 1$ on $L_{\infty}$, and inside
$B(0,3/2)\backslash\{0\}$ the limit restricts to the identity on
$E^{\perp}$.

\smallskip

Recall that by construction of the connection $d''_W$, the maps
$K(-1) \to W$, $W\to K(1)$ of the monad (3.1) are holomorphic
with respect to this connection in $B(0,2)$.  Thus, in this ball
the connection $\rho A$ on $E$ is induced by Hermitian projection
from the monad
$$
M'': \qquad 0 \longrightarrow K(-1)
\buildrel{A}\over\longrightarrow \buildrel{d''_W} \over W
\buildrel{B}\over\longrightarrow  K(1) \longrightarrow 0\quad
,\eqno(5.3)
$$
where  the notation indicates that $W$ is equipped
with the connection $d''_W$.  Since $g \cdot d''_W= d_W$ in
$B(0,3/2)$,  it follows that $gA(z) \colon K(-1) \to W$ and
$B(z)g^{-1} \colon W \to K(1)$ are holomorphic (with respect to
the standard $\db$-operators).  Indeed the connection $\rho A$ on
$E$ is precisely the ``pull-back" of the induced hermitian
connection on the cohomology of the monad
$$
M':\qquad 0
\longrightarrow K(-1)
\buildrel{gA}\over{\quad\longrightarrow\quad} \buildrel{d_W}
\over W \buildrel{Bg^{-1}}\over{\quad\longrightarrow\quad} K(1)
\longrightarrow 0\quad ,\eqno(5.4)
$$
meaning that if $E(M'')$,
$E(M')$ are the cohomologies of the respective monads, then
$$
E(M'') =ker\,B/im\,A \owns e+AK(-1) \mapsto ge+gAK(-1) \in
ker\,Bg^{-1}/im\,gA=E(M')
$$
defines a smooth map such that
$(\db_{\rho A}=)\; \db''= g\circ \db'\circ g^{-1}$ and $(\d_{\rho
A}=)\; \d''= {g^*}^{-1}\circ \d' \circ g^*$, where the hermitian
structures on each of the bundles is that which is induced by
hermitian projection from the flat metric on $W$. Choosing a
unitary isomorphism between $E(M'')$ and $E(M')$ gives a
hermitian connection on $E(M'')$  of the form $g_1\cdot \rho A$
for some complex automorphism $g_1$. Since the sequences
$\{A(z)\},\, \{B(z)\}$  converge to limits which are
non-degenerate away from the origin and the sequences $\{g\},
\,\{g^{-1}\}$ are uniformly bounded on compact subsets of the
complement of the origin, the same holds for $\{g_1^{}\}$ and
$\{g_1^{-1}\}$.

The fact that the automorphisms $g$ in the sequence are
converging to the identity on $E^{\perp}$ away from the origin
implies that the holomorphic maps $gA$ and $Bg^{-1}$ are
converging to the standard degenerate forms $A_1 :=\bmatrix{z^0&
z^1 &0}^{T}$, $B_1 := \bmatrix{z^1 & -z^0 &0}$ respectively. A
priori this convergence is on compact subsets of
$B(0,3/2)\backslash\{0\}$, but since these maps are holomorphic
it follows from Cauchy's Theorem~(Removable Singularities) that
the convergence is in fact throughout all of $B(0,3/2)$.

\medskip

The next step is to further modify the automorphisms $g$ so that
the holomorphic maps $\t A := gA$ and $\t B := Bg^{-1}$ have a
more convenient form.

\bigskip
\noindent{\bf Proposition~5.5.} \quad {\sl  Let
$D:=\{(z^0,z^1)\in \m C^2\mid |z^0|<1, |z^1|< 1\}$ be the unit
polydisk and let $A_1 := \bmatrix{z_0{\bf 1}&z^1{\bf 1}&0}^T \in
Hom(K,W)$, $B_1;=\bmatrix{z_1{\bf 1}&-z^0{\bf 1}&0}\in Hom(W,K)$.
If  $\epsilon := \sup_{D}[\,|\t A-A_1|+|\t B-B_1|\,]$ is
sufficiently small  there exists a holomorphic matrix $h$ on $D$
together with a constant matrix $\beta \in End \,K$ with
$\sup_D|h-1| + |\beta-1| \le C_0\epsilon$ such that $h\t
A=A_1+A_0$ and $\t Bh^{-1}=\beta(B_1+B_0)$ for some constant
$A_0$ and $B_0$.}
\medskip
\noindent {\bf Proof:} \quad  The
proof boils down to an application of the Implicit Function
Theorem, for which a little preparation is required.
\smallskip
Let ${\cal H}$ be the set of functions holomorphic in a
neighbourhood of $\bar D$ and let $\b{\cal H}$ be its completion
under the norm $||*||$ given by  $||f|| := \sup_D|f|$, so
$\b{\cal H}$ is a Banach algebra. If $f \in {\cal H}$, there is a
unique function $f_1 \in {\cal H}$ such that
$f(z^0,z^1)-f(z^0,0)=z^1f_1(z^0,z^1)$.  The maximum of $|f_1|$ on
$D$ is attained at a point $(a,b)$ where $|a|=1$ or $|b|=1$, but
if $|a|=1$, then since $f_1(a,z^1)$ is holomorphic in $z^1$ this
function attains its maximum on the unit disk at a point where
$|z^1| =1$, so $|f_1|$ always attains its maximum at a point
where $|z^1|=1$.  Hence $||f_1||=||z^1f_1||\le 2||f||$.  If
$\{f_i\} \in {\cal H}$ is a Cauchy sequence converging to $f\in
\b{\cal H}$, then $\{(f_i)_1\}$ is Cauchy and converges to $f_1$.
It follows that the map $\b{\cal H}\owns f \mapsto f_1 \in
\b{\cal H}$ is a $C^1$ map. Similarly, for a function $f=f(z^0)$
alone, the decomposition $f=f(0)+z^0f_0(z^0)$ gives a
differentiable function  $\b{\cal H}\owns f \mapsto f_0 \in
\b{\cal H}$.

Now let ${\cal A}$, ${\cal B}$ respectively be the spaces of
linear maps $K\to W,\,W\to K$ with coefficients in ${\cal H}$,
and let ${\cal G}$ be the space of endomorphisms of $W$ with
coefficients in ${\cal H}$ which, with respect to the
decomposition $W=[K \oplus K\oplus R]^T$ have the block form
$$
h=\bmatrix{h_{00}&h_{01}&h_{02}\cr
h_{10}&h_{11} & h_{12}\cr
h_{20}&h_{21}&h_{22}}= \bmatrix{ h_{00}\hfill& h_{01} &
h_{02}\hfill\cr
h_{10,0}(z^0)+z^1h_{10,1}(z^0)& h_{11} &
h_{12}(z^0)\cr
h_{20}(z^0)\hfill & h_{21} &0\hfill}\;. \eqno(5.6)
$$

Let ${\cal C}$ be the space of triples $(A,\beta_0,\beta_2)\in
Hom(K,W)\oplus Hom(K,K)\oplus Hom(R,K)$    with coefficients in
${\cal H}$ such that $A(0,0)=0$, $\beta_2(0,0)=0$, and if
$A=\bmatrix{z^0\alpha_{0,0}(z^0) +
z^1\alpha_{0,1}&z^0\alpha_{1,0}(z^0) +
z^1\alpha_{1,1}&\alpha_2}^T$ then
$\beta_0(z^0,z^1)=z^0\alpha_{1,0}+z^1\beta_{0,1}$ for some
$\beta_{0,1}$ satisfying $\beta_{0,1}(0,0)=-\alpha_{0,0}(0)$.
Define $P \colon {\cal A}\times {\cal B} \to C$ by
$$
P(A,B):=
\big(A-A(0)-A_1,z^0\alpha_{1,0} +z^1({\bf
1}_K-\alpha_{0,0}(0)+\beta_{0,1}-\beta_{0,1}(0,0)),
\beta_2-\beta_2(0)\big)
$$
for $A=\bmatrix{z^0\alpha_{0,0}(z^0)
+z^1\alpha_{0,1}&z^0\alpha_{1,0}(z^0) +
z^1\alpha_{1,1}&\alpha_2}^T$ and
$B=\bmatrix{\beta_0&\beta_1&\beta_2}$.

Now define a map $F$ from a neighbourhood of $(A_1,B_1,0)\in
{\cal A}\times {\cal B}\times {\cal G}$ to ${\cal C}$ by
$$
F(A,B, h) := P((1+ h)A,B(1+ h)^{-1})\;.
$$
Then
$F(A_1,B_1,0)=(0,0,0)$  and the partial derivative of $F$ in the
${\cal G}$-direction at $(A_1,B_1,0)$ is given by
$$
{\cal G}
\owns  h \mapsto \bigg(\bmatrix{z^0h_{00}+z^1h_{01}\cr
z^0h_{10}+z^1h_{11}\cr
z^0h_{20}+z^1h_{21}},
-(z^1h_{00}-z^0h_{10}), -(z^1h_{02}-z^0h_{12})\bigg) \in {\cal C}
\eqno(5.7)
$$
after a short calculation to check on the second
component.

If $(A,\beta_0,\beta_1) \in {\cal C}$ has the  form
$$
\bigg(
\bmatrix{z^0\alpha_{0,0}+z^1\alpha_{0,1}\cr
z^0\alpha_{1,0}+z^1\alpha_{1,1} \cr
z^0\alpha_{2,0}+z^1\alpha_{2,1}}, z^0\alpha_{1,0}+z^1\beta_{0,1},
z^0\beta_{2,0}+z^1\beta_{2,1}\bigg)
$$
for some $\beta_{0,1}$
such that $\beta_{0,1}(0,0)=-\alpha_{0,0}(0)$ then there is a
unique solution $h$ of the form (5.6) mapped by (5.7) onto
$(A,\beta_0,\beta_1)$, namely
$$
h=\bmatrix{
\alpha_{0,0}-z^1\beta_{0,11}\hfill &
\alpha_{0,1}+z^0\beta_{0,11}\hfill & -\beta_{2,1}\hfill \cr
\alpha_{1,0}+z^1(\alpha_{0,00}+\beta_{0,10}) &
\alpha_{1,1}-z^0(\alpha_{0,00}+\beta_{0,10}) &
\hphantom{-}\beta_{2,0}\hfill \cr
\alpha_{2,0}\hfill &
\alpha_{2,1} \hfill & \hphantom{-}0\hfill}
$$
where
$\alpha_{0,0}=\alpha_{0,0}(0)+z^0\alpha_{0,00}(z^0)$ and
$\beta_{0,1}=
\beta_{0,1}(0,0)+z^0\beta_{0,10}(z^0)+z^1\beta_{0,11}$.

Applying the Implicit Function Theorem, if $\epsilon := ||\t
A-A_1||+||\t B-B_1||$ is sufficiently small, there is a
holomorphic matrix $h$ with $||h|| \le Const.\epsilon$ such that
$P((1+h)\t A,\t B(1+h)^{-1})=(0,0,0)$.  This means that $(1+h)\t
A-A_1$ is constant, that the third component of $\t B(1+h)^{-1}$
is constant, and that the first component of $\t B(1+h)^{-1}$ has
the form $\beta_0(z^0)+z^1\beta_{01}$ where $\beta_{01}$ is
constant.
\medskip
Write $(1+h)\t A=\bmatrix{z^0{\bf 1}_K+\tilde
a_0&z^1{\bf 1}_K+\tilde a_1&\tilde a_2}^T$ and  $\t B(1+h)^{-1} =
\bmatrix{\beta_0+z^1\beta_{01}&\beta_1&\tilde b_2}$ where
$\beta_{01},b_2$ are constant and $\beta_0$ is a function of
$z^0$ alone.  Since $1+h$ is close to $1$ and the first component
of $\t B$ is close to $z^1{\bf 1}_K$, the matrix $\beta_{01}$ is
close to the identity;  similarly, $\tilde a_0, \tilde a_1,
\tilde a_2, \tilde b_2$ are close to $0$.

The monad equation $\t B\t A=0$ implies that
$$
(\beta_0+z^1\beta_{01})(z^0{\bf 1}_K+\tilde a_0)+ \beta_1(z^1{\bf
1}_K+\tilde a_1)+\tilde b_2\tilde a_2=0\,;
$$
i.e.,
$$
[\beta_0(z^0{\bf 1}_K+\tilde a_0)+\beta_{1}(z^0,0)\tilde
a_1+\tilde b_2\tilde a_2]+z^1[\beta_{01}(z^0{\bf 1}_K+\tilde
a_0)+\beta_1+\beta_{1,1}\tilde a_1]=0\;. \eqno(5.8)
$$
where $\beta_1$ has been decomposed into
$\beta_1(z^0,0)+z^1\beta_{1,1}$. The first term in brackets on
the left of (5.8) is a function of $z^0$ alone, so it follows
that each of the two terms must vanish separately. The vanishing
of the second term implies
$$
[\beta_{01}(z^0{\bf 1}_K+\tilde
a_0)+\beta_1(z^0,0)]+ \beta_{1,1}(z^1{\bf 1}_K+\tilde a_1) = 0\;.
\eqno(5.9)
$$

To proceed requires a small lemma:

\bigskip
\noindent{\bf Lemma~5.10.} {\sl Suppose $c\in GL(k,\m
C)$, $v_0 \in \m C^k$ is a row vector, and $v_1$ is a row vector
of holomorphic functions of the single complex variable $z$  in
an open set $\Omega \subset \m C$. If $v_0 + v_1(z{\bf
1}+c)\equiv 0$ in $\Omega$ and  all eigenvalues of $c$ lie within
$\Omega$  then it follows $v_0=0\equiv v_1$.}

\medskip
\noindent {\bf Proof:} \quad Choose $g\in GL(k,\m C)$
such that $g^{-1}cg := u$ is upper triangular.  Then
$v_0g+v_1g(z{\bf 1}+u) \equiv 0$.  The diagonal entries in $u$
are the eigenvalues of $c$ and by working successively along the
the components of $v_0g$, $v_1g$, the vanishing of each component
in turn follows by taking $z$ to be minus the appropriate
eigenvalue.  \qquad \qed
\medskip
Since $\tilde a_1$ is small, so
too are its eigenvalues and therefore they can be assumed to be
well within the unit disk. Hence the  lemma can be applied to
(5.9) (holding $z^0$ fixed and taking $z=z^1$) to give
$\beta_{1,1}=0$ and $\beta_1(z^0,0)=-\beta_{01}(z^0{\bf
1}_K+\tilde a_0)$, so $\beta_1$ is the linear function
$-\beta_{01}(z^0{\bf 1}_K+\tilde a_0)$. Substituting this into
the vanishing of the first term in (5.8) gives the equation  $
\beta_0(z^0{\bf 1}_K+\tilde a_0)-\beta_{01}(z^0{\bf 1}_K+\tilde
a_0)\tilde a_1+\tilde b_2\tilde a_2 =0 $, or
$$
(\beta_0-\beta_{01}\tilde a_1)(z^0{\bf 1}_K+\tilde a_0)+\tilde
b_2\tilde a_2+ \beta_{01}(\tilde a_1 \tilde a_0 -\tilde a_0
\tilde a_1)=0\;.
$$
Since $\tilde a_0$ is small Lemma~5.10 can be
applied again to  obtain $\beta_0=\beta_{01}\tilde a_1$ and
$\tilde b_2\tilde a_2+\beta_{01}(\tilde a_1 \tilde a_0 -\tilde
a_0 \tilde a_1)=0$.  Thus
$$
\t Bh^{-1}=\beta_{01}\bmatrix{{\bf
1}_K+\tilde a_1&-(z^0{\bf 1}_K+\tilde a_0)&\beta_{01}^{-1}\tilde
b_2^{}}\;,
$$
which completes the proof of the proposition.
\qquad \qed

\bigskip
Consider now the new monad  $ M_1': \quad 0
\longrightarrow K(-1)
\buildrel{\;(1+h)gA\;\;}\over\longrightarrow \buildrel{d_W}\over
W\buildrel{\;\; Bg^{-1}(1+h)^{-1}\;}\over\longrightarrow  K(1)
\longrightarrow 0 $ constructed from (5.4) using Proposition~5.5.
As in the case of the monads $M''$ and $M'$, the automorphism
$(1+h)g$ induces an isomorphism between the cohomologies $E(M'')$
and $E(M'_1)$ such that the connection on $E(M'')$ induced by
hermitian projection from the flat metric on $W$ and the
connection $d''_W$ (i.e., $d_{\rho A}$ in $B(0,2)$) is the
``pull-back" of the connection $d_1'$ on $E(M_1')$ induced by the
flat metric on $W$ and the connection $d_W$.

Since the constant terms in the new holomorphic maps $\t A$, $\t
B$  are close to $0$, $M_1'$ can only degenerate near the origin.
But by construction, it is non-degenerate in the unit polydisk,
and therefore it defines a non-degenerate monad on all of $\m
P_2$. Thus the cohomology $E(M_1')$  is a holomorphic $r$-bundle
on $\m P_2$ with $c_2=k$ and which is trivial on $L_{\infty}$.
By Theorem~0.1 of [B4] there is a smooth automorphism $g_1'$ of
$E(M_1')$ such that $d_2' := g_1'\cdot d_1'$ is anti-self-dual
(with respect to the flat metric on $\m C^2$) with $g_1'^*g_1'
=1$ on $L_{\infty}$.

Composing $g_1'$ with the map on cohomology induced by $(1+\chi
h)g$ for some appropriate cutoff function $\chi$ and then fixing
a unitary isomorphism between $E(M'')$ and $E(M_1')$, it follows
that there is a smooth automorphism $g_2$ of $E$ such that in
$B(0,1/2)$, $g_2\cdot \rho A$ is the restriction  of an instanton
on $S^4$, $\h F(g_2\cdot \rho A)$ is everywhere uniformly bounded
in $C^0$, and $g_2^*g_2^{}=1$ on $L_{\infty}$.  Since
$\lap(\log|g_2|^2) \le 2(|\h F(\rho A)|+|\h F(g_2\cdot \rho A)|)
\le Const.$ with a similar equation for $g_2^{-1}$, it follows
from the Maximum Principle that $g_2$ and its inverse remain
uniformly bounded in $C^0$ as the sequence concentrates.  By
Lemma~2.1 the sequence of instantons defined by the monads $M_1'$
converges everywhere in $B(0,1/2)$ except at the origin, and on
$B(0,1/2)\backslash \{0\}$  the limit $g_{2,\infty}$ gives an
isomorphism between the holomorphic bundle defined by the limit
of the connections $\rho A$ and the limiting instanton.  Part 4
of that lemma therefore implies that the sequence of instantons
defined by the monads $M_1'$ is bubbling all of its curvature at
the origin and degenerates nowhere else.  Moreover, if the main
theorem can be proved for the sequence of instantons
$\{E(M_1')\}$, then because the maps $g_2$ and their inverses
remain uniformly bounded, the lemma implies that the same
convergence property is enjoyed by the sequence of connections
$\{\rho A\}$.

The proof of Statement 3 of Theorem~1.3 will therefore be
essentially complete if it can be shown that a sequence of
instantons on $\m C^2$ which is bubbling all of its charge at the
origin can be made to converge by blowing up and pulling back.
This is achieved in the next section.

\bigskip
\bigskip
\bigskip

\noindent{\bf 6. \quad Convergence of monads.}
\bigskip

As stated in \S3, every instanton on $S^4$ can be described by a
monad of the form
$$
\m{P}_2:\qquad 0 \longrightarrow K(-1)
\buildrel{A}\over\longrightarrow W
\buildrel{B}\over\longrightarrow K(1) \longrightarrow 0\quad
,\eqno(6.1)
$$
where $A= \bmatrix{z^0{\bf 1}_K+a_0&z^1{\bf
1}_K+a_1&a_2}^T$,  $B=\bmatrix{z^1{\bf 1}_K+a_1&-(z^0{\bf
1}_K+a_0)&b_2}$ with the constants $a_i,b_j$ satisfying the monad
condition $a_1a_0-a_0a_1+b_2a_2=0$ and the instanton condition  $
a_0^{}a_0^*-a_0^*a_0^{}+a_1^{}a_1^*-
a_1^*a_1^{}+b_2^{}b_2^*-a_2^*a_2^{}=0
$. The monad is the restriction to $\m P_2$ of a monad on $\m
P_3$ of the form (6.1), $\m P_3$ being the twistor space for
$S^4$ ([AHS]).
\medskip
There is an analogous description of
(self-dual) instantons on $\m {CP}_2$ in  terms of monads on the
blowup  $\t{\m P}_2$ of $\m P_2$ at a point.  The twistor space
for $\m{CP}_2$ is the flag manifold $\m F = \{(Z,W)\in \m
P_2^{}\times \m P^*_2 \mid Z\cdot W=0\}$ where $Z=(Z^0,Z^1,Z^2),
\;W=(W_0,W_1,W_2)$ are homogeneous coordinates.  The hypersurface
$W_2=0$ is isomorphic the blowup $\t{\m P}_2$ of $\m P_2$ at
$(0,0,1)$, realised by the projection on the first factor, and is
also isomorphic to the  the Hirzebruch surface $H_1$, realised as
a  $\m P_1$ bundle over $\m P_1$ by projection onto second
factor. On this hypersurface, the equation $Z\cdot W=0$ implies
that $(Z^0,Z^1)=\lambda(-W_1,W_0)$ for some section $\lambda$ of
${\cal O}(1,-1)$, this section defining the exceptional line.

Holomorphic bundles on $H_1$ were studied in [B2]. With $x:=
c_1({\cal O}(1,0))$ and $y := c_1({\cal O}(0,1))$ being  first
Chern classes of the the pull-backs of the Hopf line bundles,  a
holomorphic $r$-bundle $E$ on $\t{\m P}_2$ with $c_1(E)= l(x-y)$
and $c_2(E)=kxy$ which is trivial on $L_{\infty}=\{Z^2=0=W_2\}$
is given by a  monad of the form
$$
\t{\m P}_2:\qquad 0
\longrightarrow K_1(0,-1) \buildrel{(a,b)}\over\longrightarrow N
\oplus K_2(1,-1) \buildrel{(c,d)}\over\longrightarrow\ K_3(1,0)
\longrightarrow  0\quad ,\eqno(6.2)
$$
where $K_1,\,K_2,\,K_3$
and $N$ are complex vector spaces of  dimension $k + {1\over
2}l(l-1),\, k + {1\over 2}l(l+1),\,k +  {1\over 2}l(l-1)$ and $r
+ k + {1\over 2}l(l-3)$ respectively.  The results of [B4], or
alternatively the direct results of King [Ki]  when $c_1(E)=0$,
imply that each such monad can be extended to a  unitary monad on
$\m F$, these unitary monads being described in  [B1].  That is,
the monad (6.2) is isomorphic to a monad of the same  form on $\m
F$ for which the map $(c,d)$ is
$$
(c(Z),d(W)) = (a(\bar
Z)^*,-b(\bar W)^*)\;.\eqno(6.3)
$$
Again, the failure of the map
$(a(W),b(Z))$ to be injective at some point $(Z,W)$ corresponds
precisely to the singularity of the corresponding self-dual
connection at the point $[\bar Z \times W] \in\m {CP}_2$; for
details, see [B1].

\medskip

Given a monad on $\t {\m P}_2$ of the form (6.2), a monad on $\m
P_2$ of the form (6.1) is constructed by taking direct images,
the  latter monad being non-singular iff $b(0,0,1)$ is an
isomorphism  (and the original is non-singular).  Of more
interest here is the  construction of a monad on $\t{\m P}_2$
from one on ${\m P}_2$.  Given a monad on $\m P_2$ of the form
(6.1) with cohomology $E$,  construct a monad on $\t{\m P}_2$ of
the form (6.2) corresponding  to $\pi^*E$ as follows: set $K_1 :=
K =: K_3$, $K_2 := Im\,A_2$ where $A_2 := A(0,0,1)$ and $N :=
K_2^{\perp}$.  Let $\Pi_{\bullet}$ denote hermitian projection,
and set  $a(W):=\Pi_{N}A(-W_1,W_0,0)$, $b(Z) := \Pi_{K_2}A(Z)$,
$c(Z)=B(Z)\mid_N$ and $d(W) := B(-W_1,W_0,0)\mid_{K_2}$.

The unitary condition (6.3) is not satisfied in general, but it
is possible to write down, in a relatively explicit way, the
isomorphic monad on $\t{\m P}_2$ which does satisfy the unitary
condition and possesses the same holomorphic trivialisation on
$L_{\infty}$.  Given the monad (6.1), let $a(W), b(Z)$ be as
constructed in the preceding paragraph.  The new maps $a'(W),
b'(Z), c'(Z),d'(W)$ are then given by
$$
\bmatrix{a'\cr
b'}=\bmatrix{\phi_N & 0 \cr
\lambda\phi_{KN} & \phi_K}
\bmatrix{a\cr
b}g_1^{-1}\;,\qquad \bmatrix{c'&d'}=g_2\bmatrix{c
&d} \bmatrix{\phi_N & 0 \cr
\lambda\phi_{KN} & \phi_K}^{-1}
$$
for some $\phi_N\in GL(N)$, $\phi_K\in GL(K_2)$, $\phi_{KN}\in
Hom(N,K_2)$, $g_1\in GL(K_1)$ and $g_2\in GL(K_1)$.  Setting
$\phi :=\bmatrix{\phi_N&0\cr\phi_{KN}&\phi_K} \in GL(W)$, this is
$$
\bmatrix{a'(W)\cr
b'(Z)}=\bmatrix{\Pi_N\phi A(W_1,-W_0,0) \cr
\Pi_{K_2} \phi A(Z^0,Z^1,Z^2)}g_1^{-1}\;, \quad  \bmatrix{c'(Z)
&d'(W)}= g_2\bmatrix{B(Z)\phi^{-1}\mid_N &
B(W_1,-W_0,0)\phi^{-1}\mid_{K_2}} \;.
$$

The requirement that the cohomologies of the two monads have the
same trivialisations on $L_{\infty}$ is equivalent to the
condition that $\phi$  be the identity on the subspace $R \subset
W$, which implies that in terms of the decomposition $W = K\oplus
K \oplus R$ it should have the form
$$
\phi=\bmatrix{\phi_{00}&\phi_{01}&0\cr\phi_{10}&\phi_{11}&0\cr
\phi_{20}&\phi_{21}&1}\;. \eqno(6.4)
$$

Restricting to $L_{\infty}$,  the unitary condition (6.3) on the
new monad is
$$
\bmatrix{\Pi_N\phi A(W_1,-W_0,0)g_1^{-1}\cr
\Pi_{K_2} \phi A(Z^0,Z^1,0)g_1^{-1}}=\bmatrix{(g_2B(\bar W_0,\bar
W_1,0)\phi^{-1}\mid_N)^* \cr
-(g_2B(\bar Z^1,-\bar
Z^0,0)\phi^{-1}\mid_{K_2})^*}\;;
$$
that is, $\phi^*\phi
A(Z^0,Z^1,0)=B(-\bar Z^1,\bar Z^0,0)^*g_2^*g_1^{}$, implying
$\phi^*\phi=\bmatrix{g_2^*g_1^{}&0&*\cr
0 & g_2^*g_1^{} &*\cr0 &
0 & *}$.  Using (6.4), it follows that there is a unitary matrix
$U \in GL(K\oplus K)$ together with a positive definite matrix $P
\in GL(K)$ such that
$$
\phi=\bmatrix{U_{00}P^{1/2}&U_{01}P^{1/2}
&0\cr
U_{10}P^{1/2}&U_{11}P^{1/2}&0\cr
0&0&1}\,\qquad {\rm and
}\qquad g_2^*g_1^{}=P\;.
$$

The condition that $\phi$ should preserve the subspace $K_2=Im\,
A_2$ for $A_2 := A(0,0,1)$ means that $\phi A_2=A_2\chi$ for some
$\chi \in GL(K)$;  $\chi$ must be $\psi_0^{-1}A_2^*\phi A_2$ for
$\psi_0 := A_2^*A_2^{}$.   The unitary condition (6.3) implies
one more constraint on the matrices $\phi, g_1, g_2$, namely that
if $a'(0,0,W_2) := c'(0,0,\bar W_2)^*$ and $d'(0,0,W_2) := -
b'(0,0,\bar W_2)^*$ then
$c'(0,0,1)a'(0,0,1)+d'(0,0,1)b'(0,0,1)=
c'(1,0,0)a'(1,0,0)+d'(1,0,0)b'(1,0,0)$.
A straightforward calculation shows that this condition is
equivalent to
$$
a_0^{}P^{-1}a_0^*+a_1^{}P^{-1}a_1^*+
a_0^{}Q^{-1}a_0^*+a_1^{}Q^{-1}a_1^*+b_2^{}b_2^*
=P^{-1}+P^{-1}QP^{-1} \eqno(6.5)
$$
for $Q := \chi^*\psi_0\chi$.

Note that the equation $\phi A_2 = A_2\chi$ implies $\phi^*\phi
A_2=\phi^* A_2\chi$, so $\chi^*A_2^*\phi = A_2^*\phi^*\phi$.
Hence  $ \chi^*\psi_0\chi=\chi^*A_2^*\phi A_2= A_2^*\phi^*\phi
A_2$, or
$$
(Q=)\quad
\chi^*(a_0^*a_0^{}+a_1^*a_1^{}+a_2^*a_2^{})\chi
=a_0^*Pa_0^{}+a_1^*Pa_1^{}+a_2^*a_2^{}\;. \eqno(6.6)
$$
(In fact,
the condition $\phi A_2=A_2\chi$ and the form of $\phi$ imply
$a_2=a_2\chi$, so the terms involving $a_2$ can be dropped from
this equation).

\medskip

Suppose now that the monad (6.1) is a typical element in a
sequence of such (satisfying the instanton condition) such that
$a_0,a_1,a_2,b_2$ are converging to $0$.  The monads on $\t{\m
P}_2$ constructed above (satisfying the unitary condition) then
correspond to a sequence of self-dual instantons on $\m{CP}_2$.
By weak compactness,  after gauge transformations,  there is a
subsequence of these connections together with a  finite set of
points where the curvatures concentrate, converging smoothly off
this set of points.  In terms of the monads, this means that
there are unitary automorphisms $U_N, U_{K}, U_{K_2}$ such that a
subsequence of $(U_N^{} a'U_K^{-1},U_{K_2}^{}b'U_K^{-1})$
converges to define a monad on $\t{\m P}_2$ which is degenerate
only at a finite set of points.  Thus if $U_W
:=\bmatrix{U_N&0\cr0&U_{K_2}}$, then a subsequence of the maps
$U_W^{}\phi A(Z)U_K^{-1}$ converges  to define a  monad which is
degenerate at only finitely many points in $\t{\m P}_2$. By
passing to a subsequence it can be assumed that the unitary maps
$U, U_W, U_K$ converge.

Restricting now to a point on $L_{\infty}$ which is {\it not\/}
a point of degeneration, it follows that the automorphisms
$P^{1/2}g_1^{-1}$ converge in $GL(K)$, so in fact there is no
degeneration on $L_{\infty}$.  Restricting next to a point in $\m
P_2\backslash\{L_{\infty}\cup\{(0,0,1)\}\}$ which is not a point
of degeneration, it follows that the maps $P^{1/2}a_0^{}g_1^{-1},
P^{1/2}a_1^{}g_1^{-1}$ and $a_2^{}g_1^{-1}$ converge. A point
$(z^0,z^1) \in \m C^2\backslash\{0\}$ is a point of degeneration
if there is a non-zero vector $v\in K$ such that
$lim\,(P^{1/2}g_1^{-1}z^0+P^{1/2}a_0g_1^{-1})v,
lim\,(P^{1/2}g_1^{-1}z^1+P^{1/2}a_1g_1^{-1})v$ and $lim\,
a_2g_1^{-1}v$ all vanish.  But since
$\det(P^{1/2}g_1^{-1}z^0+P^{1/2}a_0g_1^{-1})=
\det(P^{1/2}g_1^{-1})\det(z^0{\bf
1}_K+g_1^{}a_0g_1^{-1})=\det(P^{1/2}g_1^{-1})\det(z^0{\bf
1}_K+a_0)$ which converges to $\det(P^{1/2}g_1^{-1})(z^0)^k$,
with a similar statement for $a_1$, it follows that there are in
fact no points of degeneration  in $\m C^2\backslash\{0\}$.  Thus
if the sequence of monads on $\t{\m P}_2$ is to degenerate, it
can do so only at points of the exceptional line.

\smallskip
Now let $\tilde a_0 :=
P^{1/2}a_0(\chi^*\psi_0\chi)^{-1/2}, \tilde a_1 :=
P^{1/2}a_1(\chi^*\psi_0\chi)^{-1/2}, \tilde a_2 :=
a_2(\chi^*\psi_0\chi)^{-1/2}$ so  $\tilde a_0^*\tilde a_0^{}+
\tilde a_1^*\tilde a_1^{}+ \tilde a_2^*\tilde a_2^{} =
(\chi^*\psi_0\chi)^{-1/2}[a_0^*P a_0^{}+  a_1^*P a_1^{}+a_2^*
a_2^{}](\chi^*\psi_0\chi)^{-1/2}={\bf 1}_K$ by (6.6).  Passing to
a subsequence, it can be assumed that the maps $\tilde a_i$ all
converge, and the limit must satisfy the same identity.

Suppose that the sequence of monads on $\t{\m P}_2$ is bubbling
{\it all\/} of its charge at a single point on the exceptional
line.  For the sake of argument, let this point be the point with
$Z=(0,0,1)$, $W=(0,1,0)$.  Then if $A_0 := A(1,0,0)$ the
conditions that $a'(0,1,0), b'(0,0,1) \to 0$ are
$$
\phi
A_0g_1^{-1}-A_2^{}\psi_0^{-1}A_2^*\phi A_0g_1^{-1} \to 0, \quad
\psi_0^{-1/2}A_2^*\phi A_2g_1^{-1} \to 0\;. \eqno(6.7)
$$
Since
$A_2^*\phi={\chi^*}^{-1}A_2^*(\phi^*\phi)$ it follows that
$A_2^*\phi A_0g_1^{-1}
={\chi^*}^{-1}A_2^*(\phi^*\phi)A_0g_1^{-1}=
{\chi^*}^{-1}a_0^*Pg_1^{-1}$,
and  since  $U^{-1}\bmatrix{a_0\cr
a_1}= \bmatrix{P^{1/2}a_0\cr
P^{1/2}a_1}\chi^{-1}$, multiplying the first equation in (6.7)
through by the inverse of $\bmatrix{U&0\cr
0&1}$ implies
$$
\bmatrix{P^{1/2}g_1^{-1}\cr
0\cr
0} - \bmatrix{P^{1/2}a_0\cr
P^{1/2}a_1\cr
a_2}(\chi^*\psi_0\chi)^{-1}\bmatrix{a_0^*P&a_1^*P&a_2^*}
\bmatrix{1\cr0\cr0}g_1^{-1}
\to 0\,,
$$
or
$$
\bmatrix{P^{1/2}g_1^{-1}-\tilde a_0\tilde
a_0^*P^{1/2}g_1^{-1} \cr
-\tilde a_1\tilde
a_0^*P^{1/2}g_1^{-1}\cr
-\tilde a_2\tilde a_0^*P^{1/2}g_1^{-1}}
\to 0\;.
$$
Finally, since $P^{1/2}g_1^{-1}$ is converging to a
non-singular automorphism of $K$, multiplying  on the right by
the inverse of $P^{1/2}g_1^{-1}$ shows that the first equation in
(6.7) is equivalent to
$$
\bmatrix{{\bf 1}_K-\tilde a_0^{}\tilde
a_0^* \cr
-\tilde a_1^{}\tilde a_0^*\cr
-\tilde a_2^{}\tilde
a_0^*} \to 0\;;
$$
(the second condition is equivalent to
$\psi_0^{1/2}\chi g_1^{-1} \to 0$, but this is not required
here). Thus $\tilde a_0$ converges to an element of $U(K)$ and
$\tilde a_1, \tilde a_2$ converge to $0$.

In general, if the sequence  is bubbling all of its charge at a
point on the exceptional line where $[W_0,W_1,0]=[w_0,w_1,0]$ for
$|w_0|^2+|w_1|^2=1$, then it follows that $\tilde a_2 \to 0$ and
$(\tilde a_0,\tilde a_1) \to (\tilde aw_1,-\tilde a w_0)$ for
some $\tilde a \in U(K)$.

Recall that $\tilde a_i= P^{1/2}a_i(\chi^*\psi_0\chi)^{-1/2}$ for
$i=0,1$.  Thus if each term in the original sequence on $\m P_2$
is such that $\det a_0=0=\det a_1$ then it is {\it not\/}
possible for the corresponding sequence of monads on $\t{\m P}_2$
to bubble all of its charge at a single point.  (In particular,
if $k=1$, then the corresponding sequence on $\t{\m P}_2$ must
converge.)  If $-\lambda_i$ is an eigenvalue of $a_i$ and the
original sequence of instantons defined by (6.1) is pulled back
under the ``translation" $(Z^0,Z^1,Z^2) \mapsto (Z^0+\lambda_0
Z^2,Z^1+\lambda_1 Z^2,Z^2)$ (i.e., the isometry of the flat
metric on $\m C^2$ $(z^0,z^1)\mapsto (z^0+\lambda_0,
z^1+\lambda_1)$), then  $(a_0, a_1,a_2,b_2) \mapsto
(a_0+\lambda_0{\bf 1},a_1+\lambda_1{\bf 1}, a_2,b_2)$, and the
sequence of monads on $\t{\m P}_2$ corresponding to this
``translated" sequence could not bubble all of its charge at a
single point.

Instead of pulling back the monads on $\m P_2$ under these
``translations" and then constructing  corresponding monads on
$\t{\m P}_2$, the latter monads can be  viewed as those which are
obtained by blowing up at the points
$(-\lambda_0,-\lambda_1,1)\in \m P_2$ and constructing the
corresponding monads for each of these  blowups, {\it then\/}
switching to a fixed set of coordinates.  Iterating the entire
procedure at most $2k-1$ times, it follows that after blowing up
at a suitable (convergent) set of points, there is a convergent
subsequence of the corresponding sequence of Hermitian-Einstein
connections.

Thus, if $\omega_0$ is the standard flat metric on $\m C^2$,
there is a converging family of blowups of $\m C^2$ converging to
a blowup  $\t{\m C}{}^2 \buildrel{\pi}\over \to \m C^2$ with
exceptional divisor $\pi^{-1}(0)$, together with automorphisms
$g_{i_j}$ satisfying $g_{i_j}^*g_{i_j}^{}={\bf 1}$ on
$L_{\infty}$ and $\sup_K(|g_{i_j}^{}|+|g_{i_j}^{-1}|)\le C_K$ for
each compact set $K \subset \m C^2\backslash\{0\}$  such that
$\{g_{i_j}^{}\cdot \pi_{i_j}^*\rho A_{i_j}^{}\}$ converges
strongly on $\t{\m C}{}^2$;  moreover, for some metric of the
form $\omega_{0\alpha}$, the connections satisfy
$\omega_{0\alpha}\wedge F(g_{i_j}^{}\cdot \pi_{i_j}^*\rho
A_{i_j}^{})=0$. Note that if $\omega_{\beta}$ is any other metric
on $\t{\m C}{}^2$ which is close  to $\omega_0$ in a
neighbourhood of $L_{\infty}$ in the sense of [B3], it follows
from the main result of that reference together with the Maximum
Principle and Lemma~2.1 that the corresponding sequence of
$\omega_{\beta}$-Hermitian-Einstein connections has a strongly
convergent subsequence.

\bigskip
To complete the proof of Theorem~1.3, let $X,\,\omega,
\{A_i\}$ be as given in the hypotheses with $A_i$ converging
weakly off $S \subset X$ to define a quasi-bundle $E$.  Assume
for the moment that $E$ is stable, and replace $\omega$ by a
metric which is standard in a neighbourhood of each of the points
of $S$ and $A_i$ by the Hermitian-Einstein connection associated
to the new metric, placed in a gauge as specified by Lemma~2.1.
Using the procedure of \S\S 1--6 and passing to a subsequence if
necessary, construct a converging sequence of blowups $\t X_i$ of
$X$ consisting of at most $2C(E_{\rm top})-2C(E)-1$ individual
blowups converging to a blowup $\t X$ of $X$ so that in a
neighbourhood of each component of the exceptional divisor there
are complex automorphisms $g_i$ with $|g_i^{}|+|g_i^{-1}|$
uniformly bounded on compact subsets of the complement of this
component and with $g_i^{}\cdot\pi_i^*A_i^{}$ converging
strongly.  Replacing $g_i$ by $\exp({1\over
2}\rho\log(g_i^*g_i^{}))$ for some smooth cutoff function $\rho$
and making a gauge transformation then extends $g_i$ to the whole
of $\t X_i$ so that it is the identity on the complement of a
neighbourhood of the exceptional divisor.  Thus, using
Proposition~1.2 with Lemmas 2.1 and 2.2, the proofs of Statements
1--4 of Theorem~1.3 are complete.

The proof of Statement 6 of Theorem~1.3 is exactly the same
argument as that which followed the proof of Lemma~2.2, once it
is known that each bundle in the sequence of pulled-back bundles
is stable with respect to the same metric (or rather, its class
in $H^2(\t X, \m R)$).  If $b_1(X)$ is even, this follows
immediately from Proposition~1.2, but in general, if it were not
true, then it would be possible to construct a sequence of stable
bundles on $X$ converging weakly to a bundle which destabilises
the weak limit $E$;  (for details, see [B5, Para.~1 \S6].

\medskip

To deal with the case of weak limits which are not stable, there
is a useful technique for ``stabilising" semi-stable $r$-bundles
$E$.  Pick any $x_0\in X$  and blow up $X$ at this point.  Let
${\cal O}(1)$ be the dual of the line bundle defining the
exceptional divisor $L_0$, and in an annular region surrounding
$L_0$, identify $E$ with ${\cal O}(1)\oplus {\cal O}^{r-1}$ This
defines a bundle $\t E$ on $\t X$ restricting to  ${\cal
O}(1)\oplus {\cal O}^{r-1}$ on $L_0$ with $(\pi_*\t E)^{**}=E$
and $C(\t E)=C(E)+(r-1)/2r$.  It is not hard to show that the
isomorphism classes  of such bundles can be identified with $\m
P(E_{x_0}^*)$;  see [B5].  The following pair of results is also
proved in that reference (Prop.~4.2, Lemma~5.1):

\bigskip

\noindent{\bf Lemma~6.8.} \quad {\sl If the above operation is
applied at any $n > r^2$ distinct points in $X$, then for generic
choices of the elements of $\m P(E_{x_i}^*)$  the resulting
bundle $\t E$ on $\t X$ is $\omega_{\epsilon\alpha_0}$-stable for
all $\epsilon$ sufficiently small, where $\alpha_0=(1,1,\dots,
1)$.} \quad \qed

\bigskip
\noindent ({\bf Note: } If $r=2$, then at most $3$
points are required.)
\bigskip

\noindent{\bf Lemma~6.9.} \quad {\sl Let $\{A_i\}$ be a sequence
of irreducible Hermitian-Einstein connections converging weakly
to $A$ off $S \subset X$, and let $E_i$ be the corresponding
stable bundle.  Let $E$ be the semi-stable weak limit
corresponding to $A$, and let $\t E$ on $\t X$ be a stabilisation
of $E$ as described above with none of the blown up points lying
in $S$.   Then there are stabilisations $\t E_i$ of $E_i$ (stable
with respect to $\omega_{\epsilon\alpha_0}$ for fixed
$\epsilon>0$ sufficiently small) converging weakly to $\t E$ off
$\pi^{-1}(S)$ such that if $\t A_i$ is the corresponding
$\omega_{\epsilon\alpha_0}$-Hermitian-Einstein connection, then
for some complex automorphism $g_i$ of $E_{\rm top}$, $\t A_i =
g_i\cdot A_i$ on the complement of the exceptional divisor and
$\sup_K(|g_i^{}|+|g_i^{-1}|)$ is uniformly bounded for every
compact subset $K$ of this complement.} \quad \qed

\bigskip
\noindent  Statement 5 of Theorem~1.3 follows
immediately from these lemmas, completing the proof of that
theorem.

\bigskip
\bigskip
\bigskip

\noindent{\bf 7. \quad Bundles of rank $2$.}

\bigskip
\medskip

Suppose now that $b_1(X)$ is even, and let $\{E_i\}$ be a
sequence of stable $r$-bundles all topologically isomorphic to
$E_{\rm top}$ converging weakly to $E$ off $S \subset X$.  Let
$\{\t X_i\}$ be any sequence of blowups of $X$ converging to a
blowup $\blowup X$ of $X$, equipped with positive $\dbd$-closed
$(1,1)$-forms $\pi_i^*\omega+\rho_{\alpha}$ for suitable $\alpha
\in \m R^n$ satisfying $|\alpha| < \epsilon_0$ where $\epsilon_0$
is as in Proposition~1.2.  Although the blowups $\t X_i$ may be
changing, the usual weak compactness applies, as does
semi-continuity of cohomology (since a set with strictly
pseudo-convex boundary will have strictly pseudo-convex boundary
with respect to all nearby complex structures).  Thus there is a
finite subset $\t S \subset \t X$ such that (a subsequence of)
$\{\pi^*E_i\}$ converges weakly to a quasi-stable limit $\t E$
off $\t S$ with respect to $\omega_{\alpha}$.  Then $\det \t E =
\pi^*\det E$,  and there are non-zero homomorphisms $(\pi_*\t
E)^{**} \to E$, $E \to (\pi_*\t E)^{**}$.

For each $\epsilon\in (0,1]$ the same argument using the metric
$\pi_i^*\omega+\rho_{\epsilon\alpha}$ yields a weak limit $\t
E_{\epsilon}$. There are only finitely many possible first Chern
classes for the summands of each $\t E_{\epsilon}$, so for
sufficiently small $\epsilon$, the topological type of the
summands of $\t E_{\epsilon}$ will be independent of $\epsilon$.
Since the direct image of each summand is semi-stable (by
Proposition~1.2) and each must have the same have the same slope
with respect to $\omega$, $\pi_*\t E_{\epsilon}$ is semi-stable
for $\epsilon$ sufficiently small.

If $E$ is stable then the map $E \to (\pi_*\t E_{\epsilon})^{**}$
must be an isomorphism, so by Proposition~1.2, $\t E_{\epsilon}$
is stable with respect to $\omega_{\epsilon'\alpha}$ for all
$\epsilon' \in (0,\epsilon_0/|\alpha|)$.   Semi-continuity of
cohomology also gives a non-zero homomorphism $\t E \to \t
E_{\epsilon}$, but since the latter is $\omega_{\alpha}$-stable,
this must be an isomorphism. (Note that, regardless of the
stability of $E$, if $\t E_{\epsilon}$ is stable for some
$\epsilon$ then the weak limit is independent of $\epsilon \in
(0,\epsilon_0/|\alpha|)$.)  Thus $(\pi_*\t E)^{**}=E$, and there
are inter-twining operators $g_i$ linking the Hermitian-Einstein
connections on $E_i$ and $\pi_i^*E_i^{}$ over the complement of a
neighbourhood of the exceptional divisor which, together with
their inverses are uniformly bounded in $C^0$, and by Lemma~2.1
of [B4], $\t S\subset \pi^{-1}(S)$;  this proves the claim made
in the discussion following Proposition~1.2.  Note that
Proposition~1.2 combined with the same argument shows that the
weak limit $\t E$ is independent of the choice of suitable
$\alpha \in \m R^n$ satisfying $|\alpha|<\epsilon_0$.

\smallskip
Suppose now that $E$ is only semi-stable;  the
situation is considerably simpler if it is assumed from now on
that the rank of $E$ is $2$.  Thus $E=L_1\oplus L_2$ for some
line bundles $L_1,L_2$ on $X$ such that
$c_1(L_1)+c_1(L_2)=c_1(E_{\rm top})$ and
$\deg(L_1)=\deg(L_2)=\deg(E_{\rm top})/2$.

As noted above, if $\t E_{\epsilon}$ is stable for some
$\epsilon$, then it is in fact independent of  $\epsilon$.  If
$\t E_{\epsilon}$ is not stable, then for $\epsilon$ sufficiently
small as above it must have the form $\pi^*L_1\otimes L^* \oplus
\pi^*L_2\otimes L$ where $L$ is a line bundle on $\t X$ which is
trivial off the exceptional divisor and satisfies
$\rho_{\alpha}\cdot c_1(L)=0$.  If $\epsilon' \in
(0,\epsilon_0/|\alpha|)$ and $\t E(\epsilon')$ splits as $\t
L_1\oplus \t L_2$ then by definition of the number $\epsilon_0$
and the condition that $\omega_{\alpha}\cdot c_1(\t
L_1)=\omega_{\alpha}\cdot c_1(\t E)/2=\omega\cdot
c_1(E)/2=\omega\cdot c_1(L_1)$ it must be the case that $(\pi_*\t
L_1)^{**}=L_1$ (after renumbering if necessary).  This implies
$\t L_1\otimes \pi^*L_1{}^*$ is trivial off the exceptional
divisor and has degree $0$ with respect to $\omega_{\alpha}$.
Semi-continuity of cohomology gives a non-zero section of $\t
L_1\otimes \pi^*L_1{}^*\otimes L$, but since this line bundle has
degree $0$ it follows that it is trivial. Thus the weak limit
$\t E_{\epsilon}= \pi^*L_1\otimes L^* \oplus \pi^*L_2\otimes L$
is independent of $\epsilon \in (0,\epsilon_0/|\alpha|)$,
isomorphic to $\t E=\t E_1$ and $\pi_*\t E$ is, as before,
semi-stable.

For generic choices of $\alpha$, there are no non-trivial line
bundles $L$ on $\t X$ which are trivial off the exceptional
divisor and have degree $0$ with respect to $\omega_{\alpha}$.
For such $\alpha$, it follows that the weak limit on $\t X$ is
either $\omega_{\alpha}$-stable or is isomorphic to $\pi^*E$.
However, even in the case of non-generic $\alpha$, if $L$ is
non-trivial then $-c_1(L)^2>0$, implying $C(\t E) > C(E)$ and
that the amount of charge bubbled by the sequence on the blowup
is strictly less than that which was bubbled by the sequence on
$X$.

\smallskip
With these preliminary remarks in hand, the proof of
Theorem~1.4 can now be given;  the method is by induction on the
amount of charge bubbled by the sequence $\{E_i\}$.

\smallskip

For $E_i$ and $E=L_1\oplus L_2$ as above, pick $3$ points $T
\subset X\backslash S$ and stabilise $E$ to $E'$ on some blowup
$X'$ of $X$ centered at $T$.  Let $\alpha=\epsilon_1(1,1,1)$ be
such that $E'$ is $\omega_{\epsilon\alpha}$-stable for all
$\epsilon\in (0,1]$, and choose stabilisations $E_i'$ on $X'$
converging weakly to $E'$ off $S$ with respect to
$\omega_{\alpha}$, as prescribed by Lemma~7.3.  Now  use
Theorem~1.3 to construct a sequence of blowups $\t X_i'$ of $X'$
converging to a blowup $\t X'$ of $X'$ centered at $S$ with a
strongly converging (sub)-sequence $\{\pi_i^*E_i'\}$, stable with
respect to $\pi_i^*\omega_{\alpha}+\rho_{\beta}$ for some fixed
$\beta$ sufficiently small.  Without loss of generality it can be
assumed that $supp(\rho_{\alpha})\cap
supp(\rho_{\beta})=\emptyset$, so that
$\omega_{\epsilon\alpha,\delta\beta}
=\pi'^*\pi^*\omega+\epsilon\rho_{\alpha}+\delta\rho_{\beta}$ is a
positive $(1,1)$-form on $\t X'$ for any $(\epsilon,\delta)\in
(0,1]\times(0,1]$.

Let $\t E$ be the weak limit of the sequence $\{\pi^*E_i\}$ on
$\t X$, the latter space being that which obtained by blowing
down the components $T'$ of the exceptional divisor lying over
$T$;  let $\pi' \colon \t X' \to \t X$ be the blowing-down map.
By the discussion above, it can be assumed without loss of
generality that $\t E=\pi^*E$;  similarly, it can also be assumed
without loss of generality that the weak limit of the sequence
$\{\pi'^*\pi^*E_i\}$ on $\t X'$ is $\pi'^*\pi^*E$.

For each $i$ there is a map $E_i' \to \pi'^*E_i$ which is an
isomorphism off $T'$, so by semi-continuity of cohomology there
is a non-zero map $\t E' \to L_1\oplus L_2$, where for notational
convenience, the pull-back of $L_i$ to $\t X'$ has been denoted
by $L_i$ rather than $\pi'^*\pi^*L_i$.  If this map were an
isomorphism off $T'$ then $(\pi'_*\t E')^{**}=\pi^*(L_1\oplus
L_2)$, but this contradicts the fact that the charge on the
former bundle is strictly greater than that on $E$.  If the map
$\t E' \to L_2$ (say) is generically non-zero, then after taking
the maximal normal extension of the kernel of this map, it
follows that there is an exact sequence of the form $0 \to K_1
\to \t E' \to K_2\otimes {\cal J}\to 0$  where ${\cal J} \subset
{\cal O}_{\t X'}$ is a sheaf of ideals such that $supp({\cal
O}_{\t X'}/{\cal J})$ is a finite set,  where $K_1,\, K_2$ are
line bundles on $\t X'$ satisfying $K_1\otimes K_2 = L_1\otimes
L_2\otimes{\cal O}(-T')$, and there is a non-zero holomorphic map
$K_2 \to L_2$.  Since $(\pi_*\t E')^{**} = E'$ and
$(\pi'_*E')^{**} = E= L_1\oplus L_2$ it follows that
$\pi'_*\pi_*^{}\t E'$ is semi-stable. Then  $\omega\cdot
L_2=\omega\cdot c_1(\t E')/2 \le \omega \cdot c_1(K_2)$, and
therefore the map $(\pi'_*\pi_*^{}K_2)^{**}\to L_2$ must be an
isomorphism.  Thus $K_2=L_2\otimes K_2'$ where  $K_2'$ is trivial
off the exceptional divisor in $\t X'$;  moreover, $H^0(\t
X',K_2'^*)\not = 0$.

Since $(\pi_*\t E')^{**}=E'$, the stabilisation construction of
[B5,\S4] forces $K_2'$ to restrict to either ${\cal O}(0)$ or
${\cal O}(1)$ on each component of $T'$ with the latter occurring
at most once, and also forces $supp({\cal O}_{\t X'}/{\cal J})
\cap T' =\emptyset$.  Write $K_2'=: M_2'\otimes L$ where $L$ is a
line bundle on $\t X$ which is trivial off $\t S$ and $M_2'$
denotes the line bundle which is trivial off $T'$ and has the
same restriction to $T'$ as $K_2'$.  Then $M_1' := {\cal
O}(-T')\otimes (M_2')^*$ is trivial off $T'$, restricts to ${\cal
O}(0)$ on at most one component of $T'$ and to ${\cal O}(1)$ on
the others.  The bundle $\t E'$ is given by the exact sequence
$$
\t X': \qquad \qquad 0 \longrightarrow L_1\otimes M_1'\otimes
L^* \longrightarrow  \t E' \longrightarrow L_2\otimes M_2'\otimes
L\otimes {\cal J}\longrightarrow 0\;. \eqno(7.1)
$$
The condition
that $H^0(\t X',K_2'^*)\not = 0$ implies $L^*$ has a non-zero
section, so $\rho_{\beta}\cdot c_1(L)\le 0$.  Note that after
taking double-duals, $\pi'_*\pi_*^{}(7.1)$ gives the sequence $0
\to L_1 \to (\pi'_*\pi_*\t E')^{**} \to L_2\otimes {\cal J}'\to
0$ for some ${\cal J}'$ with $supp({\cal O}_X/{\cal J}') \subset
\pi\pi'(supp({\cal O}_{\t X'}/{\cal J}))\cup S$, but since
$(\pi'_*\pi_*\t E')^{**}=L_1\oplus L_2$ the sheaf ${\cal J}'$
must be isomorphic to ${\cal O}_X$;  hence  $supp({\cal O}_{\t
X'}/{\cal J}) \subset \t S$.

The direct image of (7.1)$\otimes{\cal O}(T')$ under $\pi'_*$
gives the exact sequence
$$
0  \longrightarrow L_1\otimes L^*
\longrightarrow (\pi'_*\t E')^{**}\longrightarrow L_2\otimes
L\otimes {\cal J}\longrightarrow 0\;. \eqno (7.2)
$$
If
$\pi'_*\t E'$ is semi-stable then $\omega_{\beta}\cdot
c_1(L_1\otimes L^*) \le (1/2)\omega_{\beta}\cdot c_1(\pi'_*\t
E')$, but since $\omega_{\beta}\cdot c_1(L_1\otimes
L^*)=(1/2)\omega_{\beta}\cdot c_1(\pi'_*\t E')$ and
$\omega_{\beta}\cdot c_1(L^*)=\rho_{\beta}\cdot c_1(L^*) \ge 0$,
it follows $\rho_{\beta}\cdot c_1(L^*)=0$.  Since $L^*$ has a
non-zero section, it must therefore be trivial, giving $0 \to L_1
\to (\pi'_*\t E')^{**}\to L_2\otimes {\cal J}\to 0$. Since
$C((\pi'_*\t E')^{**})=C(E_{top})> C(E)$, it follows that ${\cal
J} \not= {\cal O}_{\t X}$.

Desingularise the sequence $0 \to L_1\to (\pi'_*\t E')^{**} \to
L_2\otimes{\cal J}\to 0$ as in [B3,\S3] to obtain a blowup $\h X$
of $\t X$ centered at $supp({\cal O}_{\t X}/{\cal J})$ such that
$\hat\pi^*((\pi'_*\t E')^{**})$ is an extension of genuine line
bundles, where the line sub-bundle strictly destabilises
$\hat\pi^*((\pi'_*\t E')^{**})$ with respect to
$\omega_{\beta,\gamma}$ for any suitable $\gamma$.  On the other
hand, since $\t E'$ is $\omega_{\alpha,\beta}$-stable, it follows
that if $|\gamma|$ is sufficiently small, $\hat\pi^*\t E'$ on
$(\h X)'$ is $\omega_{\alpha,\beta,\gamma}$-stable, but is not
$\omega_{\epsilon\alpha,\beta,\gamma}$-semi-stable for $\epsilon$
sufficiently small.

If each blowup in the sequence $\t X_i'$ is now modified by
blowing up in such a way that the new sequence converges to $(\h
X)'$, then the usual arguments show that (for some subsequence),
the pullbacks of $E_i'$ converge strongly to $\hat\pi^*\t E'$
with respect to $\omega_{\alpha,\beta,\gamma}$.  Thus, it can be
assumed from the outset that in fact $\pi'_*\t E'$ is not
$\omega_{\beta}$-semi-stable, which implies that $\t E'$ is not
$\omega_{\epsilon\alpha,\beta}$-stable for sufficiently small
$\epsilon>0$.  In (7.1) therefore, $\rho_{\beta}\cdot L^* >0$
else $\pi'_*\t E'$ is an extension by torsion-free stable sheaves
of the same slope, which is automatically semi-stable.

\smallskip

For each small $\epsilon>0$, let $\t B_{\epsilon}'$ be the weak
limit of the sequence $\pi^*E_i'$ using the metric
$\omega_{\epsilon\alpha,\beta}$.  Since there are only finitely
many different splitting types for $\t B_{\epsilon}'$, for
$\epsilon$ sufficiently small the splitting has topologically
constant summands, with degree independent of $\epsilon$.  But
since $c_1(\t B_{\epsilon})=1$ on each exceptional fibre of
$\pi'$, it follows that $\t B_{\epsilon}$ must be stable for
sufficiently small $\epsilon$, and semi-continuity of cohomology
implies that this weak limit is in fact independent of such
$\epsilon$:  $\t B_{\epsilon}'=:\t B'$.  Semi-continuity of
cohomology also gives a non-zero homomorphism $\t E' \to \t B'$
which must have non-zero kernel at each point (since the two
bundles have the same determinants), and by Proposition~1.2,
$(\pi'\t B')^{**} =: \t B$ is $\omega_{\beta}$-semi-stable.

Since $\t B'$ is $\omega_{\epsilon\alpha,\beta}$-stable for all
$\epsilon$ sufficiently small, the composition $L_1\otimes
M_1'\otimes L^* \to \t E' \to \t B'$ must be identically zero, so
there is a non-zero homomorphism $L_2\otimes M_2'\otimes L \to\t
B'$. The maximal normal extension of $L_2\otimes M_2'\otimes L
\to\t B'$ in $\t B'$ must be of the form $L_2\otimes N_2' \otimes
P$ where $N_2'$ is trivial off $T'$ and $P$ is trivial off $\t
S$, for the fact that $\t B'$ is
$\omega_{\delta\alpha,\beta}$-stable implies (by Proposition~1.2)
that $\pi_*^{}\pi'_*\t B'$ is semi-stable.   Thus  there is a
commutative diagram with exact rows of the form
$$
\matrix{0&\longrightarrow& L_1\otimes M_1'\otimes
L^*&\longrightarrow& \t E'&\longrightarrow & L_2\otimes
M_2'\otimes L\otimes{\cal J} &\longrightarrow&0\cr
&&&&\vrule
height10pt width0pt depth 5pt \downarrow &&\downarrow&&\cr
0&\longleftarrow&{\cal J}_1\otimes L_1\otimes N_1'\otimes
P^*&\longleftarrow&\t B' &\longleftarrow &L_2\otimes N_2' \otimes
P&\longleftarrow&0\cr\cr} \eqno(7.3)
$$
where $N_1' := {\cal
O}(-T')\otimes (N_2')^*$ and ${\cal J}_1$ is a sheaf of ideals
such that $supp({\cal O}_{\t X'}/{\cal J}_1)$ is a finite set.

Since there is a non-zero homomorphism $M_2'\otimes L\to N_2'
\otimes P$ it follows that $\rho_{\alpha}\cdot c_1(M_2') \le
\rho_{\alpha}\cdot c_1(N_2')$ and $\rho_{\beta}\cdot c_1(L) \le
\rho_{\beta}\cdot c_1(P)$.  Since $\t B'$ is
$\omega_{\epsilon\alpha,\beta}$-stable for small $\epsilon$,
$\rho_{\beta}\cdot c_1(P) \le 0$.  If equality were to hold, then
$\omega_{\epsilon\alpha,\beta}\cdot(c_1(\t B')-2c_1(L_2\otimes
N_2' \otimes P))=\epsilon\rho_{\alpha}\cdot (c_1(\t
B')-2c_1(N_2'))\le \epsilon\rho_{\alpha}\cdot (c_1(\t B')-2c_1(
M_2'))<0$, violating stability.  Therefore $\rho_{\beta}\cdot
c_1(P) < 0$ implying in particular $C((\pi'_*\t B')^{**})>C(E)$.

As in the case of $\t E'$, there is a non-zero holomorphic map
$\t B' \to L_1 \oplus L_2$.  Since $C((\pi'_*\t
B')^{**})>C(L_1\oplus L_2)$ and $\det \pi'_*\t B'=L_1\otimes
L_2$, the map must have non-zero kernel everywhere. If the
composition $L_2\otimes N_2' \otimes P \to \t B' \to L_1\oplus
L_2$ were $0$, there would be a non-zero map $L_1\otimes
N_1'\otimes P^* \to L_1\oplus L_2$, which is impossible since
$\omega_{\epsilon\alpha,\beta}\cdot c_1( N_1'\otimes P^*)>0$ for
small $\epsilon$.  If $L_2\otimes N_2' \otimes P \to \t B' \to
L_1$ were non-zero, then since $L_1$ and $L_2$ have the same
degree, it would follow that $L_1=L_2$.  Thus it can be supposed
without loss of generality that $L_2\otimes N_2' \otimes P \to \t
B' \to L_2$ is non-zero.

Since $\pi'_*\pi_*^{}\t B'$ is semi-stable, the image of  $\t B'
\to L_2$ is of the form $L_2\otimes K_2\otimes {\cal J}_2$ for
some (new) line bundle $K_2$ trivial off the exceptional divisor
for which $K_2^*$ has a non-zero section. Writing
$K_2=K_2'\otimes K_2''$ where $K_2'$ is trivial off $T'$ and
$K_2''$ is trivial off $\t S$ and setting $K_1'':=(K_2')^*$,
$K_1':=(K_2'')^*\otimes {\cal O}(-T')$,  there is an exact
sequence  $0\to L_1\otimes K_1'\otimes K_1''\to\t B'\to
L_2\otimes K_2'\otimes K_2''\otimes {\cal J}_2\to 0$.  Moreover,
there are  non-zero maps $N_2' \otimes P \to K_2'\otimes
K_2''\otimes {\cal J}_2$ and $ K_1'\otimes K_1'' \to N_1'\otimes
P^*\otimes {\cal J}$, so it follows $\rho_{\alpha}\cdot c_1(N_2')
\le \rho_{\alpha}\cdot c_1(K_2') \le 0$ and $\rho_{\beta}\cdot
c_1(P) \le \rho_{\beta}\cdot c_1(K_2'') \le 0$.  Stability of $\t
B'$ with respect to $\omega_{\epsilon\alpha,\beta}$ for all small
$\epsilon$ implies $\epsilon\rho_{\alpha}\cdot
c_1(K_2')+\rho_{\beta}\cdot
c_1(K_2'')>(1/2)(\epsilon\rho_{\alpha}+\rho_{\beta})\cdot c_1(\t
B')=-(3/2)\epsilon\epsilon_0$ for all small $\epsilon>0$.  Hence
$\rho_{\beta}\cdot c_1(K_2'')=0$ and since $(K_2'')^*$ has a
non-zero section, it follows that $K_2''$ is trivial.
Furthermore, $K_2'$ can restrict to ${\cal O}(1)$ on at most one
component of $T'$, restricting trivially to the other components.
Since $M_2'$ also satisfies this condition and there is a
non-zero map $M_2' \to N_2'$, on each component of $T'$ the
latter can restrict only to ${\cal O}(1-a)$ for some $a\ge 0$,
with $a \ge 1$ on at least $2$ components.  The existence of the
non-zero map $N_2' \to K_2'$ then implies exactly the same
behaviour for $N_2'$;  moreover, if $K_2'$ does restrict to
${\cal O}(1)$ on some component of $T'$, then the same is true
for $N_2'$ and $M_2'$ on that component and $K_2'=N_2'=M_2'$.

\smallskip
The following commutative diagram with exact rows and
columns summarises the situation for $\t B'$;  in this diagram,
the sheaf $Q$ is the pushout:
$$
\matrix{&&0&&0&&0&&\cr
&&\uparrow&&\vrule height10pt width0pt depth 5pt
\uparrow&&\uparrow&&\cr
0&\longrightarrow& L_2\otimes N_2'
\otimes P&\longrightarrow& L_2\otimes K_2'\otimes {\cal
J}_2&\longrightarrow & Q &\longrightarrow &0\cr
&&{}_{||}
&&\vrule height10pt width0pt depth 5pt \uparrow&&\uparrow&&\cr
0&\longrightarrow& L_2\otimes N_2' \otimes P&\longrightarrow& \t
B'&\longrightarrow & L_1\otimes N_1'\otimes P^*\otimes {\cal J}_1
&\longrightarrow&0\cr
&&\uparrow&&\vrule height10pt width0pt
depth 5pt \uparrow&&\uparrow&&\cr
0&\longrightarrow
&0&\longrightarrow& L_1\otimes K_1'&= & L_1\otimes K_1'
&\longrightarrow&0\cr
&&\uparrow&&\vrule height10pt width0pt
depth 5pt \uparrow&&\uparrow&&\cr
&&0&&0&&0&&\cr}
$$

The skyscraper sheaves ${\cal O}_{\t X'}/{\cal J}_1$, ${\cal
O}_{\t X'}/{\cal J}_2$ must both be supported in $T' \cup \t S$,
and if ${\cal J}'_i $ is the sheaf of ideals which equals ${\cal
J}_i$ near $\t S$ and ${\cal O}$ elsewhere, the Riemann-Roch
formula gives $dim({\cal O}_{\t X'}/{\cal J}'_2)=dim({\cal O}_{\t
X'}/{\cal J'}_1)-c_1(P)^2 >0$;  thus, $supp({\cal O}_{\t
X'}/{\cal J}_2) \cap \t S \not= \emptyset$.

As was done for $\t E'$ above, desingularise the sequence $0\to
L_1\otimes K_1'\to\t B'\to L_2\otimes K_2'\otimes {\cal J}_2\to
0$ as in [B3] by blowing up $\t X'$ along the points of
$supp({\cal O}_{\t X'}/{\cal J}_2)\cap \t S$ to express the
pull-back of $\t B'$ to this blowup ($\t X'_1 \buildrel \pi_1
\over\rightarrow \t X'$ say) as an extension by genuine line
bundles (except perhaps near $T'$).   For any suitable $\gamma$,
the new line sub-bundle strictly destabilises $\pi'_*\t B'$ with
respect to $\omega_{\beta,\epsilon\gamma}$ for all $\epsilon >0$,
which implies that $\pi_1^*\t B'$  is
$\omega_{\epsilon\alpha,\beta,\gamma}$-stable for some small
$\epsilon>0$, but is not
$\omega_{\epsilon'\alpha,\beta,\gamma}$-stable for some
$\epsilon' < \epsilon$.

Fix $\epsilon>0$ and $\gamma$ with $|\gamma|$ and $\epsilon$
chosen sufficiently small that  $\omega_{\alpha,\beta,\gamma}$ is
suitable and so that for $c \in H^2(\t X_1',\m Z)$, the
constraints
$$
\omega \cdot c=0,\qquad |\rho_{\beta}\cdot c| \le
2(\epsilon|\alpha|+|\gamma|)\sqrt{C(E_{\rm top})+1},\qquad
-c\cdot c \le 4(C(E_{\rm top})+1) \eqno(7.4)
$$
imply
$\rho_{\beta}\cdot c=0$.   Now repeat the same construction as
before to obtain weak limits on $\t X_1'$: using the metric
$\omega_{\alpha,\beta,\gamma}$, the pull-backs of $E_i'$ must
converge strongly to $\pi_1^*\t E'$; with respect to
$\omega_{\epsilon\alpha,\beta,\gamma}$  for any sufficiently
small $\epsilon$ they converge to some new stable bundle $\t
B_1'$.

There are two cases to consider; namely, $\pi_{1*}\t B_1'$ on $\t
X'$ is or is not $\omega_{\epsilon\alpha,\beta}$-semi-stable.  In
the former case, since $\t B'$ is
$\omega_{\epsilon\alpha,\beta}$-stable and there is a non-zero
holomorphic map $(\pi_{1*}\t B_1')^{**}\to \t B'$, it follows
$(\pi_{1*}\t B_1')^{**}=\t B'$.  Since $\t B_1'$ is
$\omega_{\epsilon\alpha,\beta,\gamma}$-stable but $\pi_1^*\t B'$
is not, it follows in this case that $\t B_1'$ must be
non-trivial on the exceptional fibres of $\pi_1$ and therefore
$C(\t B_1')>C(\t B')$ ([B5, Prop. 2.9]).

Assuming as usual that the amount of charge bubbled by the
pull-backs of $E_i$ to $\t X_1'$ is the same as the original, the
same process as above can be repeated until such time as the
direct image of $\t B_i'$ onto $\t X_{i-1}'$ is not semi-stable
and the situation is that of the second possibility.

If $\pi_{1*}\t B_1'$ is not
$\omega_{\epsilon\alpha,\beta}$-semi-stable, then for some
$\delta \in (0,1)$  there is a line bundle $R$ on $\t X_1'$ with
a sheaf inclusion $R \to \t B_1'$ with torsion-free quotient such
that $\omega_{\epsilon\alpha,\beta,\delta\gamma}\cdot( c_1(\t
B_1')-2c_1(R))=0$. The same arguments as above lead to the
conclusion that $R$ must be of the form $L_1$ or $L_2$ tensored
with a line bundle trivial off the exceptional divisor, and that
$\rho_{\epsilon\alpha,\beta, \delta\gamma}\cdot( c_1(\t
B_1')-2c_1(R))=0$.  Stability  at $\delta =1$ implies
$\rho_{\gamma}\cdot( c_1(\t B_1')-2c_1(R))>0$, so
$\rho_{\epsilon\alpha,\beta}\cdot( c_1(\t B_1')-2c_1(R))<0$.
Stability with respect to $\omega_{\epsilon\alpha,\beta,\gamma}$
for all $\epsilon>0$ implies $\rho_{\beta,\gamma}\cdot  ( c_1(\t
B_1')-2c_1(R)) \ge 0$.  Since $|\rho_{\beta}\cdot( c_1(\t
B_1')-2c_1(R))| \le (\epsilon|\alpha|+|\gamma|) ||c_1(\t
B_1')-2c_1(R)||$, by (7.4) and choice of $\epsilon$ and
$|\gamma|$ the bundle $R$ must satisfy $\rho_{\beta}\cdot
c_1(R)=0$, so to summarise,
$$
\eqalign{ \omega\cdot ( c_1(\t
B_1')-2c_1(R))=0,\qquad &\rho_{\alpha}\cdot ( c_1(\t
B_1')-2c_1(R))<0\cr
\rho_{\beta}\cdot ( c_1(\t B_1')-2c_1(R))=0,
\qquad &\rho_{\gamma} \cdot ( c_1(\t B_1')-2c_1(R))>0\;.
}\eqno(7.6)
$$

As in the case of $\t B'$, there is a non-zero homomorphism $\t
B_1' \to L_1\oplus L_2$. The composition $R \to \t B_1' \to
L_1\oplus L_2$ cannot be identically zero, else there is a
non-zero map $R^*\otimes det(\t B_1') \to L_1\oplus L_2$ which
contradicts the fact that
$\omega_{\epsilon\alpha,\beta,\gamma}\cdot (c_1(L_i)-c_1(R^*
\otimes det(\t B_1'))) <0$ for $\epsilon$ sufficiently small.
Hence on each component of $T'$, $R$ must restrict to ${\cal
O}(a)$ for $a \ge 0$, and since $\rho_{\alpha}\cdot c_1(R)>
(1/2)\rho_{\alpha}\cdot c_1(\t B_1')$ it follows that $R$ can
only restrict to ${\cal O}(1)$ or ${\cal O}(0)$ with the former
occurring at most once.

Semi-continuity of cohomology gives a non-zero holomorphic map
$(\pi_{1*}\t B_1')^{**} \to \t B'$, but since the former is
unstable and the latter is stable, this cannot be an isomorphism.
The induced map $(\pi_{1*}R)^{**} \to \t B'$ must be identically
$0$ since $\omega_{\epsilon\alpha,\beta}\cdot(c_1(\t
B_1)-2c_1(R))<0$, so there is an induced map $(\pi_{1*}R^*\otimes
det(\t B_1'))^{**}\to \t B'$. If the composition of this map with
$\t B' \to L_1\otimes N_1'\otimes P^*\otimes {\cal J}_1$ were
identically $0$, there would be  an induced non-zero map
$(\pi_{1*}R^*\otimes det(\t B_1'))^{**}\to L_2\otimes N_2'\otimes
P$, and hence a non-zero section of $R\otimes L_1^* \otimes
N_1'^* \otimes P$.  This however is not possible since
$\omega_{\epsilon\alpha,\beta}\cdot c_1(R\otimes L_1^* \otimes
N_1'^* \otimes P)=\epsilon\rho_{\alpha}\cdot c_1(R\otimes
N_1'^*)+\rho_{\beta}\cdot c_1(P)$ is negative for sufficiently
small $\epsilon$.  Thus $R$ has the form $L_2\otimes
S_2'\otimes\hat L$ for some $S_2'$ on $X'$ trivial off $T'$ and
some $\hat L$ on $\t X_1'$ which is trivial off the exceptional
fibres of $\pi_1$, using here the facts that $\rho_{\beta}\cdot
c_1(R)=0$ and $\beta$ is assumed generic.  Moreover $S_2'$
restricts to ${\cal O}(1)$ on at most one component of $T'$ and
to ${\cal O}(0)$ on the others, and there is a non-zero
holomorphic map $N_2' \to S_2'$ (so as before, if $S_2'$ does
restrict to ${\cal O}(1)$ on some component of $T'$ then the same
is true of $N_2'$ and $M_2'$), and $\rho_{\gamma}\cdot \hat L
<0$.

\smallskip
By semi-continuity of cohomology once more, there is a
non-zero holomorphic map $\pi_1^*\t E' \to \t B_1'$.  Pull back
the sequence (7.1) to $\t X_1'$ and quotient out
$\pi_1^*(L_2\otimes M_2'\otimes L\otimes {\cal J})$ by its
torsion subsheaf to obtain a new sequence $ 0 \to L_1\otimes
M_1'\otimes  L^*\otimes \hat M^* \to  \t E' \to L_2\otimes
M_2'\otimes L\otimes \hat M\otimes\hat{\cal J}\to 0$ for some
line bundle $\hat M$ trivial off the exceptional fibre of $\pi_1$
satisfying $\rho_{\gamma}\cdot c_1(\hat M) \le 0$ and some sheaf
of ideals $\hat{\cal J}\subset {\cal O}_{\t X_1'}$ such that
$supp({\cal O}_{\t X_1'}/\hat{\cal J})$ is a finite set.  Setting
$S_1' := S_2'^*\otimes {\cal O}(-T')$, there is a diagram
$$
\matrix{0&\longrightarrow& L_1\otimes M_1'\otimes L^*\otimes \hat
M^*&\longrightarrow& \t E'&\longrightarrow & L_2\otimes
M_2'\otimes L\otimes \hat M\otimes\hat {\cal J}
&\longrightarrow&0\cr
&&&& \vrule height10pt width0pt depth 5pt
\downarrow && &&\cr
0&\longleftarrow&{\cal J}_3\otimes
L_1\otimes S_1'\otimes \hat L^*&\longleftarrow&\t B'_1
&\longleftarrow &L_2\otimes S_2'\otimes\hat
L&\longleftarrow&0\cr\cr} \eqno(7.6)
$$
where ${\cal J}_3$ is
another sheaf of ideals of the same type as $\hat{\cal J}$. The
induced map $L_1\otimes M_1'\otimes L^*\otimes \hat M^* \to \t
B_1'$ must be identically $0$ since
$\omega_{\epsilon\alpha,\beta,\gamma}\cdot (c_1(\t B_1')-
2c_1(L_1\otimes M_1'\otimes L^*\otimes\hat M^*)) \le
\omega_{\epsilon\alpha,\beta,\gamma}\cdot (c_1(\t E')-
2c_1(L_1\otimes M_1'\otimes L^*) <0$ for small $\epsilon>0$.
Therefore there is an induced map $L_2\otimes M_2'\otimes
L\otimes \hat M \to \t B_1'$, and the composition of this map
with the projection $\t B_1' \to L_1\otimes S_1'\otimes \hat L^*$
must also vanish identically since $M_2'$ must restrict to ${\cal
O}(0)$ on some component of $T'$ on which $S_1'$ restricts to
${\cal O}(1)$ (implying that $\deg(L_2\otimes M_2'\otimes
L\otimes \hat M)>\deg(L_1\otimes S_1'\otimes\hat L^*)$ with
respect to some appropriate metric on $\t X'_1$). It follows that
there is  a non-zero holomorphic map $L_2\otimes M_2'\otimes
L\otimes \hat M \to L_2\otimes S_2'\otimes\hat L$, so $\hat
L\otimes \hat M^*$ has a non-zero section.  Thus
$\rho_{\gamma}\cdot \hat M \le \rho_{\gamma}\cdot \hat L <0$,
implying that $\hat M$ cannot be trivial on all the exceptional
fibres of $\pi_1$.

Since $C(\t E') \ge
-(1/4)[(c_1(L_1)-c_1(L_2))^2+2c_1(L)^2+2c_1(\hat M)^2]$, it
follows that after repeating this construction (of $\t B_1$)
sufficiently often, this process must eventually terminate and at
that point, the pull-backs of the bundles $E_i$ must have a
subsequence which bubbles a strictly smaller amount of charge
than $C(E_{\rm top})-C(E)$.  Applying the inductive hypothesis
then yields a sequence of blowups for which the pull-backs of
$E_i$ have a strongly convergent subsequence.

\medskip

A count of the maximum number of blowups required to obtain a
strongly convergent subsequence yields the crude bound $(C(E_{\rm
top})-C(E))^2+1$. This estimate can be improved using the
following argument, which is applicable for arbitrary rank $r$.

\smallskip
Suppose that for some sequence of blowups $\t X_i$
converging to $\t X$, $\pi_i^*E_i^{}$ converges strongly to a
stable bundle $\t E$ with respect to $\omega_{\alpha}$.  Let $\t
S$ be the exceptional divisor in $\t X$ and pick a set of $r^2+1$
points  $T\subset X$ disjoint from a neighbourhood of
$supp(\rho_{\alpha})$.  Let $E_0 := (\pi_*\t E)^{**}$ and
stabilise $E_0$ to a bundle $E'_0$ on a blowup $X'$ of $X$
centered at $T$. Since $\pi^*E_0=\t E$ near $T$, there is a
corresponding ``stabilisation" $\t E'$ of $\t E$, and by
construction, $\t E'$ is stable with respect to
$\omega_{\epsilon\alpha,\delta\beta}$ for all $\epsilon,\delta\in
(0,1]$. The convergence of $\pi_i^*E_i^{}$ to $\t E$ implies that
there is a corresponding sequence of stabilisations $E_i'$ of
$E_i$ such that $\pi_i^*E_i'$ converges to $\t E'$ with respect
to $\omega_{\epsilon\alpha,\delta\beta}$.

After passing to a subsequence if necessary, there is a finite
set $S' \subset X'$ such that $\{E_i'\}$ converges weakly off
$S'$ to a quasi-stable bundle  with respect to $\omega_{\beta}$.
Stability of $E_0'$ and semi-continuity of cohomology imply that
this bundle is isomorphic to $E_0'$.  Applying Lemma~2.1 on $\t
X'\backslash supp(\rho_{\alpha})$  and using the strong
convergence of the sequence $\pi^*E_i'$, it follows that $S'
\subset X'\backslash \pi(supp(\rho_{\alpha}))$;  that is, there
is no bubbling of curvature near $\pi'^{-1}(T)$ for the sequence
$\{E_i'\}$.  Furthermore, the amount of charge bubbled by this
sequence is at most $C(E_{\rm top})$.

Now apply Theorem~1.3 to obtain a  sequence $\{\hat X_i'\}$  of
blowups of $X'$ consisting of at most $2C(E_{\rm top})-1$ blowups
converging to a blowup $\hat X'$ of $X'$ for which, after passing
to a subsequence if necessary, $\hat\pi_i^*E_i'$ converges to a
stable bundle $\hat E'$ with respect to $\omega_{\beta,\gamma}$,
where $|\gamma|$ is taken to be so small that the pull-back from
$X'$ of a  bundle of charge $\le C(\t E')$ which is
$\omega_{\beta}$-stable  is $\omega_{\beta,\gamma}$-stable.  By
Proposition~1.2, $\gamma$ can be replaced by $\delta\gamma$ for
any $\delta \in (0,1)$ to obtain the same limit $\t E'$. The
exceptional divisor of $\hat X'\to X'$ is $\hat\pi^{-1}(S')$ and
$(\hat\pi_*\hat E')^{**}=E'_0$.

Semi-continuity of cohomology implies that on a sequence of joint
blowups of $\hat X_i'$ and $\t X_i$, the pull-backs of $E_i'$
converge to the pull-back of $\t E'$ with respect to
$\omega_{\alpha,\beta,\delta\gamma}$ if $\delta>0$ is small
enough, and to the pull-back of $\hat E'$ with respect to
$\omega_{\epsilon\alpha,\beta,\gamma}$ for $\epsilon>0$
sufficiently small. Semi-continuity of cohomology gives a
non-zero holomorphic map between these two pull-backs which
pushes down to a non-zero map between the double duals of the
direct images on $X'$.  But since both such push-downs are
isomorphic to the stable bundle $E_0'$, it follows that the map
is an isomorphism off the exceptional divisor, and since the two
bundles have the same determinant, it is an isomorphism
everywhere.

Now let $\hat E := (\pi'_*\hat E')^{**}$.  Then $\pi^*\hat
E=\hat\pi^*\t E$, and since the latter is
$\omega_{\alpha,\delta\gamma}$-stable for sufficiently small
$\delta$ it follows that $\hat E$ must be
$\omega_{\delta\gamma}$-stable.  Semi-continuity of cohomology
now implies that $\{\hat\pi_i^*E_i^{}\}$ has a subsequence
strongly convergent to $\hat E$ with respect to
$\omega_{\delta\gamma}$. This completes the proof of Theorem~1.4.
\quad \qed

\bigskip
\bigskip
\bigskip

\centerline{\bf REFERENCES}
\bigskip
\medskip
\tolerance=1000
\parindent=25pt
\everypar{\hangindent .55in}
\parskip=6pt
\frenchspacing
\baselineskip=12pt

\item{[AHS]~~} M. F. Atiyah, N. J. Hitchin and I. M. Singer,
{Self-duality in four dimensional Riemannian geometry}, {\it Proc.
Roy. Soc. Lond.} Ser. A 362 (1978) 425--461.
\item{[BPV]~~} W. Barth, C. Peters and A. Van de Ven, {\it Compact
Complex Surfaces}, (Springer, Berlin, Heidelberg, New York, 1984).
\item{[B1]~~} N. P. Buchdahl, ``Instantons on $\m{CP}_2$", {\it J.
Differ. Geom.} 24 (1986) 19--52.
\item{[B2]~~} N. P. Buchdahl, ``Stable 2-bundles on Hirzebruch
surfaces", {\it Math. Z.} 194 (1987) 143--152.
\item{[B3]~~} N. P. Buchdahl, ``Hermitian-Einstein connections and
stable vector bundles over compact complex surfaces", {\it Math.
Ann.} 280 (1988) 625--648.
\item{[B4]~~} N. P. Buchdahl, ``Instantons on $n\m{CP}_2$", {\it
J. Differ. Geom} 37 (1993) 669--687.
\item{[B5]~~} N. P. Buchdahl, ``Blowups and gauge fields".
Preprint (1995).
\item{[D1]~~} S. K. Donaldson, ``Instantons and geometric
invariant theory", {\it Commun. Math. Phys.} 93 (1984) 453--460.
\item{[D2]~~} S. K. Donaldson, ``Anti-self-dual Yang-Mills
connections over complex algebraic varieties and stable vector
bundles", {\it Proc. Lond. Math. Soc.} 50 (1985) 1--26.
\item{[D3]~~} S. K. Donaldson, ``Connections, cohomology and the
intersection forms of 4-manifolds", {\it J. Differ. Geom.} 24
(1986) 275--341.
\item{[D4]~~} S. K. Donaldson, ``Irrationality and the h-cobordism
conjecture ", {\it J. Differ. Geom.} 26 (1987) 141--168.
\item{[D5]~~} S. K. Donaldson, ``Polynomial invariants for smooth
four manifolds", {\it Topology} 29 (1990) 257--315.
\item{[FU]~~} D. S. Freed and K. K. Uhlenbeck, {\it Instantons and
Four-Manifolds}, MSRI Publications Vol 1.  (Springer, New York,
Berlin, Heidelberg, Tokyo, 1984).
\item{[FM1]~~} R. Friedman and J. W. Morgan, ``On the
diffeomorphism types of certain algebraic surfaces I", {\it J.
Differ. Geom.} 27 (1988) 297--369.
\item{[FM2]~~} R. Friedman and J. W.  Morgan,  ``On the
diffeomorphism types of certain algebraic surfaces II", {\it J.
Differ. Geom.} 27 (1988) 371--398.
\item{[Gau]~~} P. Gauduchon, ``Le th\'eor\`eme de l'excentricit\'e
nulle", {\it C. R. Acad. Sci. Paris} 285 (1977) 387--390.
\item{[Gie]~~} D. Gieseker, ``On the moduli of vector bundles on
an algebraic surface", {\it Ann. Math.} 106 (1977) 45--60.
\item{[GT]~~} D. Gilbarg and N. S. Trudinger, {\it Elliptic
partial differential equations of second order.} 2nd ed.
(Springer: Berlin Heidelberg New York, 1983).
\item{[Ki]~~} A. D. King, ``Instantons and holomorphic bundles on
the blown-up plane", {\it D. Phil. Thesis} Oxford (1989).
\item{[Kob]~~} S. Kobayashi, ``Curvature and stability of vector
bundles", Proc. Japan Acad. Ser. A Math. Sci. 58 (1982) 158--162.
\item{[Kot]~~} D. Kotschick, ``On manifolds homeomorphic to
$\m{CP}^2\#8\b{\m{CP}}^2$", {\it Invent. Math.} 95 (1989) 591--600.
\item{[L]~~} M. L\"ubke, ``Stability of Einstein-Hermitian vector
bundles", {\it Manuscr. Math.} 42 (1983) 245--257.
\item{[LY]~~} J. Li and S. -T. Yau, ``Hermitian-Yang-Mills
connection on non-K\"ahler manifolds". In: Mathematical Aspects of
String Theory, ed S. -T. Yau (World Scientific: Singapore, 1987).
\item{[Ma]~~} M. Maruyama, ``On a compactification of a moduli
space of stable bundles on a rational surface". In {\it Algebraic
Geometry and Commutative Algebra}, Kinokuniya, Tokyo 1988.
233--260.
\item{[Mor]~~} J. W. Morgan, ``Comparison of the Donaldson
polynomial invariants with their algebro-geometric analogues",
{\it Topology} 32 (1993) 449--488.
\item{[OSS]~~} C. Okonek, M. Schneider and H. Spindler, {\it
Vector bundles on complex projective spaces}, (Birkh\"auser,
Boston, Basel, Stuttgart, 1980).
\item{[OV]~~} C. Okonek and A. Van de Ven, ``Stable bundles and
differentiable structures on certain algebraic surfaces", {\it
Invent. Math.} 86 (1986) 357--370.
\item{[S]~~} S. Sedlacek, ``A direct method for minimizing the
Yang-Mills functional", {\it Commun. Math.  Phys.} 86 (1982)
515--528.
\item{[U1]~~} K. K. Uhlenbeck, ``Connections with $L^p$ bounds on
curvature", {\it Commun. Math. Phys.} 83 (1982) 31--42.
\item{[U2]~~} K. K. Uhlenbeck, ``Removable singularities in
Yang-Mills fields", {\it Commun. Math. Phys.} 83 (1982) 11--30.
\item{[UY]~~} K. K. Uhlenbeck and S. -T. Yau, ``On the existence
of Hermitian-Yang-Mills connections in stable vector bundles",
{\it Commun. Pure App. Math.} 39 (1986) 257--293.

\bye